\documentclass[11pt]{article}
\usepackage{amssymb,latexsym}

\newtheorem{theorem}{Theorem}[section]

\newtheorem{defi}{Definition}[section]
\newtheorem{lemma}{Lemma}[section]
\newtheorem{prop}{Proposition}[section]

\newcommand{\df}{\displaystyle\frac}

\newcommand{\hsp}{\hspace*{\parindent}}

\newcommand{\NN}{{\Bbb N}}

\newcommand{\QQ}{{\Bbb Q}}
\newcommand{\ZZ}{{\Bbb Z}}

\newcommand{\sA}{{\cal A}}

\newcommand{\sO}{{\cal O}}
\newcommand{\sP}{{\cal P}}

\newcommand{\bz}{{\bf z}}
\newcommand{\bM}{{\bf M}}
\newcommand{\bx}{{\bf x}}

\newcommand{\bZ}{{\bf Z}}
\newcommand{\bS}{{\bf S}}
\newcommand{\bE}{{\bf E}}
\newcommand{\bV}{{\bf V}}
\newcommand{\bG}{{\bf G}}
\newcommand{\bB}{{\bf B}}
\newcommand{\bL}{{\bf L}}
\newcommand{\bW}{{\bf W}}
\newcommand{\bT}{{\bf T}}
\newcommand{\bU}{{\bf U}}
\newcommand{\bQ}{{\bf Q}}
\newcommand{\by}{{\bf y}}

\newcommand{\bv}{{\bf v}}

\newcommand{\bu}{{\bf u}}

\input amssym.def
\input amssym.tex

\newcommand{\beql}[1]{\begin{equation}\label{#1}}
\newcommand{\eeq}{\end{equation}}

\renewcommand{\theequation}{\arabic{section}.\arabic{equation}}

\catcode`\@=11
\renewcommand{\section}{
        \setcounter{equation}{0}
        \@startsection {section}{1}{\z@}{-3.5ex plus -1ex minus
        -.2ex}{2.3ex plus .2ex}{\large\bf}
        }
\catcode`@=12
\makeatletter
\def\eqalignno#1{\displ@y \ta {\bf s} kip\@centering
  \halign to\displaywidth{\hfil$\@lign\displaystyle{##}$\ta {\bf s} kip\z@skip
    & $\@lign\displaystyle{{}##}$\hfil\ta {\bf s} kip\@centering
    & \llap{$\@lign##$}\ta {\bf s} kip\z@skip\crcr
    #1\crcr}}
\makeatother

\makeatletter
\def\@sect#1#2#3#4#5#6[#7]#8{\ifnum #2>\c@secnumdepth
     \def\@svsec{}\else 
     \refstepcounter{#1}\edef\@svsec{\csname the#1\endcsname.\hskip .75em }\fi
     \@tempskipa #5\relax
      \ifdim \@tempskipa>\z@ 
        \begingroup #6\relax
          \@hangfrom{\hskip #3\relax\@svsec}{\interlinepenalty \@M #8\par}%
        \endgroup
       \csname #1mark\endcsname{#7}\addcontentsline
         {toc}{#1}{\ifnum #2>\c@secnumdepth \else
                      \protect\numberline{\csname the#1\endcsname}\fi
                    #7}\else
        \def\@svsechd{#6\hskip #3\@svsec #8\csname #1mark\endcsname
                      {#7}\addcontentsline
                           {toc}{#1}{\ifnum #2>\c@secnumdepth \else
                             \protect\numberline{\csname the#1\endcsname}\fi
                       #7}}\fi
     \@xsect{#5}}
\def\@begintheorem#1#2{\it \trivlist \item[\hskip \labelsep{\bf #1\ #2.}]}
\makeatother
\begin{document}
\begin{center}
{\Large {\bf Succinct Certificates for the Solvability of
Binary Quadratic Diophantine Equations }} \\
\vspace{\baselineskip}
{\large {\em J. C. Lagarias}} \\
\vspace{.5\baselineskip}
Department of Mathematics \\
University of Michigan  \\
Ann Arbor, MI 48109-1043 \\
lagarias@umich.edu \\
\vspace{1\baselineskip}
(July 1, 2009 revision) \\
\vspace{2\baselineskip}
{\em Abstract}
\end{center}

Binary quadratic Diophantine equations are of interest from the
viewpoint of computational complexity theory.
This class of equations includes as special cases many of the known
examples of natural problems apparently occupying intermediate stages
in the $P - NP$ hierarchy, i.e., problems not known to be solvable in
polynomial time nor to be $NP$-complete, for example 
the problem of factoring integers.

Let $L(F)$ denote the length of the binary encoding of the binary
quadratic Diophantine equation $F$ given by
$ax_1^2 + bx_1 x_2 + cx_2^2 + dx_1 + ex_2 + f = 0$.
Suppose $F$ is such an equation having a nonnegative integer solution.
This paper shows that there is a proof (i.e., ``certificate'') that $F$
has such a solution which can be checked in 
$O(L(F)^5 \log L(F) \log \log L(F))$
bit operations.
A corollary of this result is that the set $\Sigma = \{F : F$ has a
nonnegative integer solution$\}$ is in the complexity class $NP$.
The result that $\Sigma$ is in $NP$ is interesting because it is known that
there are binary quadratic Diophantine equations whose smallest nonnegative
integer solution is so large that it requires time
exponential in $L(F)$ just to write this solution down in the
usual binary representation.


\setlength{\baselineskip}{1.0\baselineskip}
%
%
%

\section{Introduction}
\hsp
There has been considerable interest  in bounding the
computational complexity of various number-theoretic problems.
A particular motivation is the Rivest-Shamir-Adleman enciphering
scheme \cite{Riv78} whose resistance to cryptanalysis depends on the
apparent difficulty of factoring large integers.
Many of these number-theoretic problems can be formulated as one of
two types of problems involving Diophantine equations.
\begin{description}
\item{(i)}
Deciding whether a given Diophantine equation has an admissible solution or not.
\item{(ii)}
Exhibiting an admissible solution to such an equation when it has one.
\end{description}
Here an {\em admissible} solution denotes an integer solution which may also be
required to satisfy some side conditions characteristic of the particular problem.
The side conditions that arise are generally of the following two types:
\begin{description}
\item{(i)}
{\em Nonnegativity.}
Certain variables are required to be nonnegative.
\item{(ii)}
{\em Congruence.}
Certain variables $x_i$ are required to satisfy congruence restrictions
$x_i \equiv \alpha_i ( \bmod~\Gamma )$ where the $\alpha_i$ and
$\Gamma$ are given as input.
\end{description}
In this framework, for example, an integer $N$ is composite if and only if the
binary quadratic Diophantine equation
\beql{eq101}
(x+2)(y+2) = N
\eeq
has a solution in nonnegative integers $x, y \in \NN$.
The problem of factoring $N$ involves exhibiting nonnegative solutions to a
series of equations (\ref{eq101}), and showing that certain other equations
of the form (\ref{eq101}) are solvable.\\

There is a close relation between Diophantine equations and the theory of
computation.
The methods developed by Davis, Putman, Robinson and Matijasevic in their
solution to Hilbert's 10th problem established that for any 
recursively enumerable
set $A$ of natural numbers there is a Diophantine equation 
$P(x_1 , \ldots , x_n ) = 0$ such that
$$
x \in A \Leftrightarrow \exists ~\mbox{nonnegative $x_2 , \ldots , x_n$ 
such that $P(x, x_2 , \ldots , x_n ) = 0$}
$$
(see in particular \cite{Dav73}, \cite{DavPutRob61}, \cite{Mat70}.)
Adleman and Manders \cite{Adl76}, \cite{AdlMan76} used these methods to
establish a computational complexity theory based on the notion of recognizing
sets for which a given Diophantine equation has solutions of size bounded
by a given complexity function $\Phi$.
More precisely, they considered sets $S$ given
$$
\begin{array}{ll}
x \in S \Leftrightarrow & \exists ~\mbox{nonnegative 
$x_2 , \ldots , x_n$ with $L(x_2 ) , \ldots , L(x_n ) \leq \Phi (L(x))$} \\
~~~ \\
 & \mbox{such that $P(x, x_2 , \ldots , x_n ) = 0$}
\end{array}
$$
where $P(x_1 , \ldots , x_n ) \in \ZZ[x_1, ..., x_n]$ is a fixed Diophantine equation, $L(x)$
denotes the length of the binary integer $x$, and $\Phi (t)$ is a complexity
measure which is an increasing function of $t$.
They introduced a complexity class $D$ which is a Diophantine analogue of the
complexity class $NP$.
It consists of all relations $R \subseteq \NN^m $ specified by a Diophantine
equation $P(x_1 , \ldots , x_{m+n} ) = 0$ with $P \in \ZZ[x_1, ..., x_{m+n}]$ and
a polynomial $q(t)$ as follows:
\begin{eqnarray}~\label{eq102}
< x_1 , \ldots , x_m > \in R & \Leftrightarrow &  \exists 
~\mbox{nonnegative $y_1 , \ldots , y_n $} \nonumber  \\ 
 &  & \mbox{such that $MAX ~L(y_i ) \leq q(L(x_1 ) + \ldots + L(x_m )) $} \nonumber \\ 
 &  & \mbox{and $P(x_1 , \ldots , x_m , y_1 , \ldots , y_n ) = 0$} ~. 
\end{eqnarray}
It is immediately clear that $D \subseteq NP$.
It is an open problem whether or not $D = NP$; this is an important problem
in determining the relative computing power of Diophantine equations as
compared to that of nondeterministic Turing machines.\\

%
%
%

\subsection{Binary Quadratic Diophantine Equations}
\hsp

The class of {\em binary quadratic Diophantine equations} (BQDE's)
\beql{eq103}
ax_1^2 + bx_1 x_2 + cx_2^2 + dx_1 + ex_2 + f = 0
\eeq
is of special interest from the viewpoint of both number theory and 
complexity theory.
 From the viewpoint of number theory, this class of equations can be used to
encode the problem of factorization, the problem of solving quadratic congruences
$$
x^2 = f ~ ( \bmod ~e )
$$
(which corresponds to $x_1^2 - ex_2 - f = 0)$, of solving Pell's equation
$$
x_1^2 - dx_2^2 = 1 ~,
$$
and  problems in representation and equivalence of binary quadratic
forms.
 From the viewpoint of complexity theory this class of equations seems to
represent the borderline between tractable and intractable computational
problems.
It includes as special cases most of the known examples of natural problems
apparently occupying intermediate stages in the $P-NP$ hierarchy (i.e.,
problems not known to be in $P$ nor to be $NP$-complete) as well as
$NP$-complete problems.
For example, it known that:
\begin{description}
\item{(i)}
$S = \{p|p$ is prime$\} \in P$ (Agrawal, Kayal and Saxena \cite{AKS04}).
\item{(ii)}
$\{ \alpha , \beta , \gamma \in \NN | \exists$ nonnegative 
$x_1 , x_2$ such that
$ax_1^2 + \beta x_2 - \gamma = 0 \}$ is $NP$-complete.
(Manders and Adleman \cite{ManAdl78}).
\item{(iii)}
$\{ a, c \in \NN | \exists x_1 , x_2$ such that 
$ax_1 x_2 + x_2 = c \} \in NP \backslash co-NP$ provided $NP \neq co-NP$.
This problem is $\gamma$-complete (Adleman and Manders \cite{AdlMan77}).
\item{(iv)}
$\{a, c \in \NN | \exists x_1 , x_2$ such that
$x_1^2 - a^2 x_2^2 = c \}$ is unfaithfully random complete (Adleman and
Manders \cite{AdlMan78}, \cite{AdlMan79}).
\end{description}
All of the sets (i)--(iv) are in $D$, and hence are certainly in $NP$.
Note that the  existence of $NP$-complete sets in $D$ does not establish 
$D = NP$.
Indeed,
Adleman and Manders \cite{AdlMan76} exhibit a set in $P$ not known 
to be in $D$.\\

This paper treats the general problem of recognizing those 
binary quadratic Diophantine
equations which have nonnegative solutions, which may also be required to
satisfy given congruence side conditions.
This problem appears to be fundamentally harder computationally than any of
the special subclasses of binary quadratic Diophantine equations considered
up to now (i.e. including (i)--(iv) above) as indicated by the following 
example. 
%
%
\paragraph{Example.} ({\em Anti-Pellian Equation})
Consider the set of equations
\beql{eq104}
x^2 - dy^2 = -1
\eeq
where $d$ is given in its binary representation as input.
The equations (\ref{eq104}) are often called the {\em non-Pellian} or
{\em anti-Pellian} equations.
For the subset $d = 5^{2n+1}$, the input requires no more than $7n$ bits.
In Lagarias \cite[Appendix A]{Lag80b}, it is shown that for
$d = 5^{2n+1}$ this equation has solutions for
$n = 1,2,3, \ldots$ and that the solution $(t_1 , u_1 )$ to this equation
with minimal binary lengths $L(t_1 )$, $L(u_1 )$ is given by
\beql{eq105}
t_1 + u_1 \sqrt 5 = (2+ \sqrt 5)^{5^n} ~.
\eeq
This implies that the length of {\em any} solution $x$ to (\ref{eq104})
expressed in binary for these $d$ satisfies
$$
L(x) > \df{1}{3} 5^n ~.~~~\Box
$$ 

This example shows that there are some binary quadratic Diophantine equations
whose solutions are so large that it requires exponential space (in terms of
the length of the coefficients of the equation) to store any such solution 
in binary.
In addition this shows the set
\beql{eq106}
\Pi = \{ d \in \NN | x^2 - dy^2 = -1 ~\mbox{is solvable in integers \}}
\eeq
cannot be established to be in the complexity class $D$ by using  the 
Diophantine
equation $x^2 - dy^2 + 1 = 0$ in a relation (\ref{eq102}).
Indeed we cannot hope to show the set (\ref{eq106}) is in $NP$ by
guessing a solution $x,y$ in binary and verifying it is a solution
by substitution in (\ref{eq104}), because 
this potentially requires exponential time
to check. 
These same restrictions apply to the general quadratic Diophantine equation 
(\ref{eq103}),
because the problem of recognizing the subclass (\ref{eq104}) 
is clearly in $P$.\\

%
%
%

\subsection{Main Results}
\hsp

The major result of this paper is that there exist short certificates 
of the solvability of all
binary quadratic Diophantine equations which have solutions.
By the preceding example, these certificates must sometimes verify that
solutions exist without exhibiting these solutions written in binary.
The certificates actually contain in a compact form enough information to
exactly calculate an admissible solution.\\

In order to state the main result, we need two definitions.
We consider a binary quadratic Diophantine equation (BQDE)  $F(x_1 , x_2 ) = 0$ where
$$
F(x_1 , x_2 ) = ax_1^2 + bx_1 x_2 + cx_2^2 + dx_1 + ex_2 + f,
$$
together with a congruence side condition
\begin{eqnarray*}
x_1 & \equiv & \alpha_1 ~( \bmod~\Gamma ) \\
x_2 & \equiv & \alpha_2 ~( \bmod~ \Gamma )
\end{eqnarray*}
with  $0 \leq \alpha_1 , \alpha_2 < \Gamma$.
We define the {\em length}
$L(F)$ of the input to be
\beql{eq107}
L(F) := L(a) + L(b) + L(c) +L(d) + L(e) + L(f) + 3L( \Gamma ) ~,
\eeq
in which  
\beql{107a}
L(a) := 2 + \log_2(|a|+1)
\eeq
 is a measure of the binary length of 
an integer $a$, allowing one extra bit for its sign.
The other definition concerns the measurement of the running time of
a program.
We shall measure running time  in terms of {\em elementary operations},
which consists of a Boolean operation on a bit or pair of bits,
and an input or shift of a single bit.
Our main result is the following, concerning  {\em nonnegative solutions} to the system
above. 

%
%
%

\begin{theorem}~\label{th11}
 Let $F(x_1 , x_2 ) = 0$ be a binary quadratic Diophantine equation, where
$$
F(x_1 , x_2 ) = ax_1^2 + bx_1 x_2 + cx_2^2 + dx_1 + ex_2 + f ~,
$$
which has an  integer solution $(x_1 , x_2 )$ satisfying the congruence condition
\begin{eqnarray*}
x_1 & \equiv & \alpha_1 ~( \bmod ~\Gamma ) \\
x_2 & \equiv & \alpha_2 ~( \bmod ~ \Gamma )~.
\end{eqnarray*}
and the nonnegativity condition
$$
x_1 \ge0, ~~x_2 \ge 0.
$$
Then there exists a certificate showing that $F(x_1 , x_2 ) = 0$ has
such an admissible solution which requires $O(L(F)^5 \log L(F) \log \log L(F))$
elementary operations to verify.
\end{theorem}

This result gives certificates imposing two side conditions on the solutions: a {\em congruence
condition} and a {\em nonnegativity condition}. 
This theorem  is formulated 
to impose  a
nonnegativity  side condition,  in order for it to provide a result compatible  
with the framework of Hilbert's 10-th problem, and also with the Diophantine
complexity theory of Adleman and Manders \cite{AdlMan77}, \cite{AdlMan78}
which requires nonnegative variables. \\

An immediate consequence of the form of the certificates produced by
Theorem \ref{th11} is the following result. 

%
%
%

\begin{theorem}~\label{th12}
 The following sets $\Sigma_i$ are all in $NP$.
\begin{center}
\begin{tabular}{rl}
(i) ~$\Sigma_1 =$ & $ \{a,b,c,d,e,f, \alpha_1 , \alpha_2 , 
\Gamma \in \ZZ |~ \exists$ {\em nonnegative integers} $x_1$, \\
 & $x_2$ {\em with} $x_1 \equiv \alpha_1 ( \bmod ~ \Gamma ), 
x_2 \equiv \alpha_2 ( \bmod ~ \Gamma )$ {\em and} \\
 & $ax_1^2 + bx_1 x_2 + cx_2^2 + dx_1 + ex_2 + f = 0 \}$ \\
~~~ \\
(ii) ~$\Sigma_2 =$ & $\{a,b,c,d,e,f \in \ZZ | ~\exists $ 
{\em nonnegative integers} $x_1 , x_2$ {\em with} \\
 & $ax_1^2 + bx_1 x_2 + cx_2^2 + dx_1 + ex_2 + f = 0\} $ \\
~~~ \\
(iii) ~$\Sigma_3 =$ & $\{a,b,c,d,e,f \in \ZZ | ~\exists$ 
{\em integers} $x_1 , x_2$ {\em with} \\
 & $ax_1^2 + bx_1 x_2 + cs_2^2 + dx_1 + ex_2 + f = 0 \} ~.$
\end{tabular}
\end{center}
\end{theorem}

Since one can tell in polynomial time whether or not a 
binary quadratic Diophantine
equation (\ref{eq103}) is of the special form
\beql{eq108}
ax_1^2 + ex_2 + f = 0
\eeq
and since the set of equations of the form (\ref{eq108}) 
which have nonnegative integral
solutions is $NP$-complete \cite{ManAdl78}, we conclude that:
{\em testing membership
in the sets $\Sigma_1$ and $\Sigma_2$ in Theorem~1.2 are each  $NP$-complete.}
We are unable to decide whether or not any of the sets 
$\Sigma_i$ above are in the
Diophantine complexity class $D$.\\

The proofs of Theorems \ref{th11} and \ref{th12} use the theory of binary quadratic forms 
in the form
developed by Gauss in Disquisitiones Arithmeticae \cite{Gau66}.
A  treatment of this theory  can be found in Buell \cite{Bu89}.
The certificates are based on Gauss' operation of {\em composition}
of forms, and crucially use an idea  of Shanks \cite{Sha72}, 
which he called the ``infrastructure'' of quadratic forms. Shanks did
not give detailed proofs of his ``infrastructure" method, but 
the  ``infrastructure" method was  put on a rigorous footing by Lenstra \cite{Le82} in 1980,
 in the framework of quadratic number fields. 
 This paper gives an alternate justification 
 of the infrastructure method in the framework of composition of forms (Lemma~\ref{le63}).
 The proof of Theorem~1.1 is outlined in Section~2, and the details appear
in the following sections. \\

In addition to the main theorem, we show that whenever two
integer binary quadratic forms are equivalent there exist succinct
certificates verifying this equivalence (Theorem~\ref{th71}).\\

%
%
%

\subsection{Related Work}
\hsp
  
  We add some remarks on related work.
The following result is  a direct consequence of Lagarias \cite[Theorem 1.1]{Lag80b}. 
%
%
%

\begin{theorem}~\label{th13}
The set
$$
\Pi_{AP}= \{ d | ~\exists ~\mbox{integers $x,y$ with $x^2 - dy^2 = -1 \}$}
$$
is in NP $\cap$ co-NP.
\end{theorem}

The problem of characterizing
the set $\Pi_{AP}$ of solvable anti-Pellian equations 
has been extensively studied in algebraic number theory,
see Narkiewicz \cite[pp. 124--126]{Nar74},  and Redei \cite{Red53}.
For recent work, see Williams \cite{Wil02},
Jacobson and Williams \cite{JW09} and Fouvry and Kl\"{u}ners \cite{FK09}.\\

In another direction Gurari and Ibarra \cite{GurIba79} consider a class of
Diophantine equations containing (\ref{eq108}) as a special case,
but not containing (\ref{eq103}).
They show that the subclass of such equations having nonnegative 
solutions is in $D$, i.e.,
when such equations have a nonnegative solution, they have one that is 
small enough
to serve as a certificate.
Their $NP$-completeness result then follows from \cite{ManAdl78}.\\

%
%
%

\subsection{Retrospective: Smale's Problem 5}
\hsp

%
%
%


These certificates given in Theorem~\ref{th11} are  also relevant to 
 Problem 5 of the  mathematical problems for the twenty-first century
 formulated   by  S. Smale \cite[p. 275]{Sm00},  which concerns height bounds for
 solutions to Diophantine equations.  \\
 
 \noindent{\bf  Smale's Problem 5.} 
 {\em  Can  one decide 
 if a single Diophantine equation $f(x, y)=0$ with
 in two variables of exact total  degree $d$, and of genus at least one, 
 has an integer solution, in time $O(2^{s^c})$ for some universal constant $c$,
 where $c$ is a universal constant? That is, can the problem be decided in exponential time?}\\
 
  In this problem
 $$
 f(u,v) = \sum_{\alpha_1 + \alpha_2 \le d} a_{\alpha_1, \alpha_2} u^{\alpha_1} v^{\alpha_2} \in \ZZ[u,v],
 $$ 
 with $f(u,v)$ having some nonvanishing  term of total degree $\alpha_1+\alpha_2= d$, and  
 $$
 s= s(f):= \sum_{\alpha_1 + \alpha_2 \le d}  \max (\log |a_{\alpha_1, \alpha_2}|, 1).
 $$
 is a measure of the ``height " of $f$.\\

Here we have the following result, which affirmatively answers Problem 5 for
polynomials of total degree $d \le 2$. 

%
%
%

\begin{theorem}~\label{th14} 
Let $F(x_1 , x_2 ) = 0$ be a binary quadratic Diophantine equation,given by
$$
F(x_1 , x_2 ) = ax_1^2 + bx_1 x_2 + cx_2^2 + dx_1 + ex_2 + f ~,
$$
given with coefficients encoded in binary. 
Then there exists a deterministic exponential time algorithm which decides whether or
not  $F(x,y)=0$ has   
an integer solution.
This algorithm uses at most   $O\left(2^{c_1L(F)} \right) $
elementary operations, where $c_1$ is an absolute constant. 
\end{theorem}

Theorem~\ref{th14} shows we may take the universal constant $c=1$ in 
Smale's problem 5,  when the input
is  restricted to bivariate polynomials  $f(u,v)$ of total degree $d \le 2$.
The anti-Pellian example above shows that to write down a minimal solution in
binary may sometimes require at least $\Omega \left( 2^{c_2 L(F)} \right)$ bits. \\

Theorem~\ref{th14}  is proved using a complexity analysis of
the  classical algorithmic approach  of Gauss for finding solutions of binary quadratic
 Diophantine equations. It  is given in Section \ref{sec56}, and  does not require
 use of the succinct certificates found in Theorem ~\ref{th11}.\\

A more general version of Theorem \ref{th14}, which works for testing
 for admissible solutions to  a BQDE that also satisfy  a congruence side condition and 
 a positivity side condition, 
 can be proved directly from Theorem~\ref{th11}.   This  (inefficient) algorithm
 used  sequentially tests all possible
candidate certificates produced in Theorem~\ref{th11}. If none of them work,
then the system has no solution.
 The number of 
candidate certificates can be shown to be $O \left( 2^{c_2L((F)^3 }\right)$ and leads to
a running time bound $O\left(2^{c_3 L(F)^3}\right).$  We omit details.\\



  In connection with his Problem 5, 
  Smale \cite[p. 276]{Sm00} also  put forward  the following hypothesis, for dealing with
curves of genus one or larger:\\
 
\noindent  {\bf  Height bound hypothesis:} {\em If the curve $f$, of positive genus, has
 any   integer solution, then it has a solution $(a, b)$ satisfying the
 estimate: $\log \max(|a|, |b|)$ is polynomially bounded by $s(f)$.}\\
 
 
 The truth of this (unproved) hypothesis would solve Smale's 
 Problem 5 affirmatively in case the genus is $1$ or larger. 
 Such a height bound does not hold for genus $0$ curves,
 as shown by  the  anti-Pell equation  example above (given in Lagarias \cite{Lag80b}).
 Effective bounds on size of solutions are known  for a large class of curves covered by
 Runge's method, see for example Walsh \cite{Wa92}, and these may imply
 the height bound in these cases. There are also very large effective bounds  
 for size of integer solutions for genus one curves, based
 on Baker's method, see Baker and Coates \cite{BC70} and Schmidt \cite{Sch92}.
 This method extends to curves given  using  Galois coverings,  for which see Bilu \cite{Bi97} 
Note that the  height bound hypothesis, if true, would provide certificates falling in 
 the Diophantine complexity class D,  and this   provides  renewed
motivation for  studying the complexity class D. \\

Concerning  the remaining case
of genus $0$ curves, Theorem~\ref{th14} solves problem 5 affirmatively for
quadratic Diophantine equations, which form a restricted class of 
 genus $0$ plane curves, without invoking the height bound hypothesis. 
 General bounds on the integer solutions of genus zero curves
 were given in Bilu and Poulakis{BP93} and improved in
 Poulakis (\cite{Po97}, \cite{Po98} , \cite{Po01} \cite{Po03}).
 However these bounds depend in part on Baker's method,
 and it appears that they are not strong enough to resolve Smale's Problem 5
 for all genus $0$ plane curves.  \\

%
%
%

%
%
%

\subsection{Retrospective: Infrastructure Method}
\hsp

The results of this paper were  announced in  preliminary report form
 in the 1979 FOCS conference proceedings \cite{Lag79}.
Detailed proofs were given in a 1981 Bell Laboratories technical report (\cite{Lag81});
 this work was contemporaneous with that of 
Lenstra \cite{Le82}.  Renewed motivation to publish
this work,  after a long delay, came from its relevance to Smale's problem 5. 
 The present paper is a slightly revised version of  \cite{Lag81}, which adds
 the application to Smale's problem 5, and corrects  errata. 

The basic idea of this paper exploits the Shanks infrastructure method,
which is here presented in the language of integral binary quadratic forms,
and composition of forms.  Since 1982 there
has been extensive development of the infrastructure method, mostly given in the language
of algebraic number fields and ideals. It is now a workhorse in computational number
theory, and is implemented in PARI 
We now review these developments. 

In 1982  H. W. Lenstra, Jr.  \cite{Le82} 
gave a rigorous analysis of the infrastructure method of Shanks
for the purpose of computing regulators (units in quadratic fields)
and class numbers.


In 1989 the  ``infrastructure" method was  used
by Buchmann and Williams \cite{BW89}  to give succinct
certificates  for class numbers and  approximate representation
of regulators of quadratic number fields, under the assumption
of the generalized Riemann hypothesis. Further work was done by 
Buchmann, Thiel and Williams \cite{BTW92}.
In 1994 Theil \cite{Th94} showed that verifying the value of the class
number falls in the class NP $\cap$ co-NP, assuming the truth of the 
generalized Riemann hypothesis. 
Recently the infrastructure has been
given a more precise theoretical formulation in terms of the Arakelov class group of a number field,
by  Schoof \cite{Sch08}. \\

 The infrastructure method is  well known to be 
computationally effective  in practice, 
as described in Chapter 5 of Cohen \cite{Coh93}.
It is used in computations of class numbers
and regulators of quadratic and cubic number fields
and function fields. For a recent survey on the computation of
solutions of Pell's equation, see Williams \cite{Wil02}.\\

%
%
%

\subsection{Acknowledgments}
\hsp

\paragraph{Acknowledgments.}

This work was revised while the author was supported
by NSF grants DMS-0500555 and DMS-0801029.

%
%
%

\section{Outline of the Proof}
\hsp
In this section we describe the main ideas of the proof
and establish some notational conventions.\\

For the proofs we deal  throughout with a system $F$ consisting of a binary quadratic
Diophantine equation (BQDE)  having the  Gauss {\em standard  form}:
\beql{eq201}
ax_1^2 + 2bx_1x_2 + cx_2^2 + 2dx_1 + 2ex_2 + f = 0 ~,
\eeq
with side conditions 
\beql{eq202}
x_i \equiv \alpha_i ( \bmod ~\Gamma ) , ~ i = 1,2 ,
\eeq
\beql{eq203}
x_i \geq 0 ,~ i = 1,2 ,
\eeq
The requirement that the coefficients of $x_1x_2$ and $x_1$ and $x_2$ 
be  even integers  is imposed  for compatibility with Gauss' formulation of this problem.
Any system can be brought to this form by multiplying (\ref{eq103}) by $2$.
A solution to (\ref{eq201})--(\ref{eq203}) will be called {\em admissible}.\\

We follow the approach of Gauss to finding solutions of such equations,
which is outlined in G. B. Mathews \cite[Chap. IX]{Mat61} and H. J. S. Smith \cite[Arts. 93-97]{Smi}.

Binary quadratic
Diophantine equations (\ref{eq201}) are classified as {\em definite, indefinite} or
{\em degenerate}
according to the value  of the {\em determinant}
\beql{eq204}
D = b^2 - ac
\eeq
being negative, positive and not a square, or a perfect square, respectively.
(The determinant $D$ is
just  $\frac{1}{4} \mbox{Disc}(f_Q)$, where $f_Q$ is
the polynomial $f_Q(x)= Q(x,1)= ax^2+bx+ c$.)
This classification is useful because the sets of solutions to these three
types of equations have qualitatively different behaviors.
In particular, definite and degenerate binary quadratic Diophantine equations 
with an admissible solution 
always
have admissible solutions small enough to serve directly as certificates.
(Lemma~\ref{le32}).
The crucial part of the proof concerns the case of indefinite binary
quadratic Diophantine equations. \\

Gauss \cite[Art. 216--221]{Gau66} gave a method to determine whether (\ref{eq201}) 
has any
integer solutions and if so to give a complete parametric description of all
solutions.
This method is based on his theory of integral binary quadratic forms, and in
particular on determining the equivalence or inequivalence of such forms.
Gauss' method easily extends to include the congruential side condition (\ref{eq202}),
but the positivity side condition (\ref{eq203}) adds new complications.
We follow the outline of Gauss's method in reducing the problem to that
of recognizing the {\em equivalence} of two quadratic forms.
In Section~3 we transform the problem to that of studying the 
generalized Pell equation
$$
x_1^2 - Dy_2^2 = g,
$$
with $(x_1, x_2)$   satisfying  side conditions on their signs and congruence
conditions to some modulus. 
In Section~4 (primitive) binary quadratic forms are introduced and the problem is transformed
to that of demonstrating that the reduced identity form $\tilde{I}$ of
determinant $D$ is equivalent to a particular reduced form $Q_{\rm red}$ via an
equivalence matrix $W$ having certain properties.
(See Section~4 for definitions.)
The proofs are complicated by the need to bound the size of the
least admissible solution and to keep track of the nonnegativity condition 
(\ref{eq203}) under these transformations.\\

These 
reductions have not yet addressed the main difficulty in finding
certificates of solvability, which is the possible exponentially large
size (number of binary bits) of the least admissible solution.
This difficulty is  here transformed into the possibly equally large size of
the entries of the equivalence matrix $W$ appearing
in Lemma~\ref{le43}. We show that in order to verify admissibility we
need only know (1) that $W$ gives an equivalence, that
(2)  the entries of $W$ satisfy certain congruence side 
conditions, and 
(3)  the entries of $W$ are known in floating point to sufficient accuracy to check a certain
sign condition.\\

The remainder of the proof is devoted to a detailed study of  those matrices $\bW$ that 
demonstrate  the equivalence of the reduced identity form $\tilde{I}$ and 
any particular reduced
form $Q_{\rm red}$.
In Section~5 we describe results of Gauss.
Gauss defined a notion of two reduced forms being neighbors.
If we form a graph in which the reduced forms are vertices, and edges
correspond to two reduced forms being neighbors, then Gauss
showed that this graph is a union of disjoint cycles. Furthermore  the cycle including
$\tilde{I}$ contains exactly the reduced forms equivalent to $\tilde{I}$,
which we call the {\em  principal cycle}.
These results imply that the associated equivalence matrices have 
a very special
form, which is related to the ordinary continued fraction algorithm.
However this form by itself is insufficient to produce succinct
certificates. However it is sufficient to yield an exponential
time algorithm for determining if a BQDE has an integer solution,
and we prove Theorem~\ref{th13}.\\

In Section~6 we come to the main idea of the succinct certificate proof.
This uses another set of relations between
these equivalence matrices, which comes from Gauss' operation of
{\em composition} of binary quadratic forms (Lemma~\ref{le62}).
The idea of using the action of composition of forms on the
principal cycle is due to Shanks \cite{Sha72}, who called it the
``infrastructure". The ``infrastructure"  asserts that composition is
a kind of doubling of distance on the graph of the set of forms.  Shanks did not give 
detailed proofs, but a rigorous justification of the 
``infrastructure" was later given by Lenstra \cite{Le82}, in the
language of ideals. In this paper we give an alternate rigorous justification in the
language of composition of forms, in  Lemma~\ref{le63}.
The action of composition can be combined 
with Gauss's  reduction steps 
 to find a short
sequence of composition formulae that prove the equivalence of any two given forms in the
principal cycle (Lemma~\ref{le64}). 
 In effect each 
composition causes a squaring, so that if one multiplied out all
the compositions to write down the matrix giving
the equivalence,  the resulting entries would
potentially have exponentially many digits, in terms of the input size.
However the correctness of the composition steps can be verified for 
each formula separately, avoiding this potential exponential blowup.\\

In Section~7 we further apply these composition formulae to give succinct certificates
for the equivalence of two 
(equivalent)
 indefinite binary quadratic forms.
We remark that the equivalence or inequivalence of two
definite or degenerate quadratic forms can be decided in polynomial 
time (\cite{Lag80a}.)\\

In Section~8 we complete the proof of Theorems~\ref{th11}, and then deduce
Theorem~\ref{th12} and Theorem~\ref{th13} from it. 
This is done by showing that the formulae of Lemma~\ref{le63}
can be used to check that the entries of $W$ satisfy a  {\em given side congruence condition},
and to evaluate $W$ using floating-point computations to enough accuracy
to verify a {\em sign condition on the solutions}. This is of interest if one wishes
to recognize positive solutions, as are studied in Hilbert's 10-th problem.
It requires significant extra work to establish these
extra side condition properties. 
In general it is difficult to rigorously prove results that the number of
significant digits present after a sequence of floating-point operations,
because there is the possibility of losing all significant
figures when adding two nearly equal floating-point number of opposite signs.
We are able to show this
potential cancellation effect cannot occur here, using a priori information about the 
magnitudes of the
quantities being computed at all intermediate steps of the computation.\\

Appendix~A gives bounds on the period lengths (mod $M$) of solutions
to certain second-order linear recurrences.
Appendix~B gives needed results on floating point computation,
concerning bounds on the loss of accuracy in floating-point operations.

%
%
%

\subsection{ Notations and Conventions.}

The {\em length} $L(a)$ of an integer $a$ is defined by
\beql{eq200a}
L(a) = 1 + \log_2(|a|+1).
\eeq
This measures its binary length, plus one bit for its sign. 

The {\em size} $||F||$ of the BQDE system 
(\ref{eq201})--(\ref{eq203}) is given by
\beql{eq205}
||F|| := MAX (|a|,|b|, |c|, |d|, |e|, |f|, | \Gamma |) ~.
\eeq
The size $||F||$ is related to the  length $L(F)$ in (\ref{eq107}) by
\beql{eq206}
\frac{1}{2} \log ||F|| \leq L(F) \leq 9\log ||F|| +18 ~,
\eeq
where $\log x$ denotes the natural logarithm.

We also need an analogous notion of {\em size} $|| \bf M ||$ of a matrix
${\bf M} = [m_{ij} ]$ given by
\beql{eq207}
|| {\bf M} || := MAX_{i,j} | m_{ij} | ~.
\eeq
If $M,N$ are $m \times k$ and $k \times n$ matrices, respectively,
then we have the trivial bound
\beql{eq208}
|| {\bf M N } || \leq k ~ || {\bf M} || ~ || {\bf N} || ~.
\eeq


When counting elementary operations (bit operations) we will sometimes use the function
\beql{eq210}
M(n) := n ( \log n) ( \log \log n )
\eeq
arising from the Sch\"{o}nhage-Strassen bound $O(M(n))$ for the multiplication of
two $n$ bit binary integers.
The O-symbol has the usual meaning, that $O(f(n))$ means $\leq c|f (n)|$,
where $c$ is an absolute, effectively computable positive constant, which
may differ at each occurrence of the O-symbol.

%
%
%

\section{Bounding the Size of the Least Admissible Solution}
\hsp
In this section we bound the size of the least admissible solution to 
definite
and degenerate binary quadratic Diophantine equations, and show that this
solution itself may serve as a certificate.
The difficult case is that of indefinite binary quadratic Diophantine equations.
In this case we obtain an exponential upper bound for the size of the minimal
admissible solution.

%
%
%
\subsection{Transformation to Generalized Pell Equation $y_1^2-Dy_2^2= g$.}

We start with a {\em standard form} binary quadratic Diophantine equation
\beql{eq300}
ax_1^2 + 2bx_1x_2 + cx_2^2 + 2dx_1 + 2ex_2 + f = 0 ~,
\eeq
with side conditions 
\beql{eq300a}
x_i \equiv \alpha_i ( \bmod ~\Gamma ) , ~ i = 1,2 ,
\eeq
\beql{eq300b}
x_i \geq 0 ,~ i = 1,2 ,
\eeq

We can immediately simplify (\ref{eq300}) by an invertible variable change
provided $D \neq 0$, $c \neq 0$, where $D=b^2-ac$. 
Following Mathews (\cite{Mat61}, 258--260) we introduce new variables
\beql{eq301}
\left[
\begin{array}{c}
y_1 \\
y_2
\end{array}
\right] = \left[
\begin{array}{c}
D ~~ 0 \\
b ~~ c
\end{array}
\right] \left[
\begin{array}{c}
x_1 \\
x_2
\end{array}
\right] + \left[
\begin{array}{c}
be -cd \\
e
\end{array}
\right]
\eeq
We have the identity
$$
y_1^2-Dy_2^2 = -cD(ax_1^2+ 2 bx_1x_2+ cx_2^2 + 2dx_1 + 2ex_2 + f ) + g
$$
where
\beql{eq303}
g = -c \left|
\begin{array}{c}
a~b~d \\
b~c~e \\
d~e~f
\end{array}
\right|  = -c(2bde-cd^2-ae^2).
\eeq
Thus we obtain that when (\ref{eq300}) holds, then 
\beql{eq302}
y_1^2 - Dy_2^2 = g ~,
\eeq
and we note the bound
\beql{eq304}
|g| \le  4|| F ||^4 ~.
\eeq
Inverting the system (\ref{eq301}) yields
\beql{eq305}
\left[
\begin{array}{c}
cDx_1 \\
~~~ \\
cDx_2
\end{array} \right] =
 \left[
\begin{array}{cc}
c & 0 \\
~~~ \\
- b & D \end{array}\right] 
~ \left[
\begin{array}{c}
y_1 \\
~~~ \\
y_2 \end{array} \right] + 
\left[
\begin{array}{c}
c(cd-be) \\
~~~ \\
c(ae-bd)
\end{array}
\right] ~.
\eeq
This yields  the following result. 
%
%

\begin{lemma}~\label{le31}
 Given a BQDE (\ref{eq300}) having $cD \ne 0$.
 Let $(y_1 , y_2 )$ be an integral solution to (\ref{eq302}),
 and let a modulus $\Gamma$ be given.
Then $(x_1 , x_2 )$ given by (\ref{eq305}) is a rational
solution to (\ref{eq300}).
The congruence class of $(y_1 , y_2 )$ (mod $cD \Gamma )$ determines
whether $x_1 , x_2$ is integral, and if so specifies
$(x_1 , x_2 )$ (mod $\Gamma )$.
\end{lemma}

%
%
%
\subsection{Certificates for Definite and Degenerate BQDE's}

We now bound the size of solutions to definite and degenerate binary
quadratic Diophantine equations. These bounds are strong enough that the
smallest admissible solution will directly serve as a polynomial time certificate. 

%
%

\begin{lemma}~\label{le32}
Suppose that  a given binary quadratic Diophantine equation system (\ref{eq300})-(\ref{eq302})
 is either
definite or degenerate.
If it has any admissible solutions at all, then it has an admissible solution
$(x_1 , x_2 )$ with
\beql{eq306a}
MAX ( |x_1| ,  |x_2|) \leq  8 ||F||^4.
\eeq
In particular , this solution satisfies 
\beql{eq306}
MAX ( L( x_1) ,  L (x_2) ) \leq 8 \log || F || + 8.
\eeq
\end{lemma} 

\paragraph{Proof.}
It suffices to prove (\ref{eq306a}) since (\ref{eq306}) follows on taking logarithms.
We treat several cases, of which the first is the generic case. \\

{\bf Case~1.}
$D \neq 0 , ~c \neq 0$.
Then  the change of variables (\ref{eq301}) takes integral
solutions of (\ref{eq201})  to integral solutions of
\beql{eq307}
y_1^2 - Dy_2^2 = g ~.
\eeq
In the definite case $D < 0$,  all integer solutions to (\ref{eq307}) have
$$
MAX (|y_1 |, |y_2 |) \leq \sqrt{|g|} ~.
$$
Using (\ref{eq304}) and (\ref{eq305}), we see that all integer solutions to (2.1) must have
$$
MAX (|x_1 |, |x_2 |) \leq \sqrt{6} (||F||^3 + ||F||^2) + 2 ||F||^2 
\le 6 || F ||^4 ~.
$$
This implies (\ref{eq306a}).\\

In the degenerate case $D = h^2$ is a perfect square and (\ref{eq307}) becomes
\beql{eq308}
(y_1 + hy_2 ) (y_1 - hy_2 ) = g ~.
\eeq
Hence each integer solution $(y_1 , y_2 )$ gives rise to a factorization
$g = g_1 g_2$ where
\begin{eqnarray}\label{eq309}
y_1 + hy_2 & = & g_1 \nonumber \\
y_1 - hy_2 & = & g_2 ~.
\end{eqnarray}
Solving (\ref{eq309})) for $(y_1 , y_2)$ we obtain
$$
MAX (|y_1 |, |y_2 | ) < |g|
$$
using (\ref{eq305}) again, this implies (\ref{eq306}). \\

{\bf Case~2.}
$D \neq 0 , c = 0, a \neq 0$.
Interchange $x_1$ and $x_2$, proceed as in Case~1. \\

{\bf Case~3.}
$D \neq 0, a = 0 , c = 0$.
Then (\ref{eq201}) multiplied by $b$ factorizes as
$$
2(bx_1 +e) (bx_2 + d) = 2de - bf ~.
$$
A similar argument to (\ref{eq308}), \ref{eq309}) yields
$$
MAX (|x_1 |, |x_2 |) \leq 3 || F ||^2
$$
in this case, implying (\ref{eq306a}). \\

{\bf Case~4.}
$D = 0$. Then 
$$
ax_1^2 + 2bx_1 x_2 + cx_2^2 = m (\alpha x_1 + \beta x_2 )^2
$$
where $m, \alpha , \beta$ are integers, $m \neq 0$.
Let
\beql{eq310aa}
z = \alpha x_1 + \beta x_2 ~.
\eeq
Suppose first $\alpha \beta \neq 0$.
Substituting $x_2 = \frac{z-\alpha x_1}{\beta}$ in (\ref{eq201}) we obtain
$$
m \beta z^2 + 2 ez  +2(d\beta-e \alpha) x_1+ f \beta = 0
$$
if $e \alpha - d \beta = 0$, then we obtain
$$
m \beta z^2 + 2 ez +f \beta=0.
$$
There are two solutions to this quadratic equation, both having
$$
|z| \le 3||F||^{\frac{3}{2}},
$$
and  if (\ref{eq310aa}) has any solution it has one with
$$
MAX(|x_1|, |x_2|) \le 6 ||F||^2.
$$
If $e \alpha - d \beta \neq 0$ then we obtain 
\beql{eq310}
x_1 = \frac{m \beta z^2 +2ez + f \beta}{2(e \alpha - d \beta )}
\eeq
and $x_2 = \frac{z-\alpha x_1}{\beta}$ yields
\beql{eq311}
x_2 =- \left(\frac{m \alpha z^2 + 2 dz + f \alpha}{2(e \alpha - d \beta )} \right) ~.
\eeq
Now the congruence class of $z( \bmod~ 2(e \alpha -d \beta ) \Gamma )$
determines whether the quantities $(x_1 , x_2 )$ given by (\ref{eq310}), (\ref{eq311}) are
integral, and if so specifies their congruence class
$(\bmod ~\Gamma )$.
Hence any  $(x_1 , x_2 )$ $(\bmod ~\Gamma )$
that can occur is given by some $z$ in any block of
$2(e \alpha -d \beta ) \Gamma$ consecutive values of $z$.
Next, we consider what signs of $(x_1 , x_2 )$ can occur.
The sign of $x_1$ changes when the numerator of the right side of (\ref{eq310})
changes sign, and the sign change occurs at some $z$ with
$|z| < 6 || F ||^2$.
A similar result holds for the sign of $x_2$ via (\ref{eq311}).
We conclude that  if there is an admissible solution, there will be one with
$$
|z| < 6 || F ||^2 + |2(e \alpha -d \beta ) \Gamma | < 8 || F ||^3 ~.
$$
Using (\ref{eq310}) and (\ref{eq311}) then gives (\ref{eq306a}).

Finally suppose $\alpha \beta = 0$.
If $\alpha = \beta = 0$ then (\ref{eq201}) is linear and (\ref{eq306a}) is easily
verified.
If $\alpha = 0$ and $\beta \neq 0$ then use (\ref{eq310}) and replace
(\ref{eq311}) with
$$
x_2 = \df{z}{\beta} ~.
$$
The same argument as in the case $\alpha \beta \neq 0$ now proves (\ref{eq306a}).
The case $\alpha \neq 0 , \beta = 0$ is treated similarly.
$~~~\Box$ \\

%
%
%
\subsection{Indefinite  BQDE's: Standard Form $y_1^2-Dy_2^2=g$  with side conditions}

Lemma~\ref{le32} produces polynomial size
certificates for definite and degenerate BQDE's. 
In the sequel it remains to  consider indefinite binary quadratic Diophantine
equations. \\

In the indefinite case a standard form equation has $cD \ne 0$,
hence the
variable change (\ref{eq301}) is invertible, and necessarily $D \ge 2$.
We now reduce the problem of finding admissible solutions to 
(\ref{eq201}) to that
of finding an admissible solution (suitably defined) to a generalized Pell equation
$y_1^2 - Dy_2^2 = g$. 
%
%

\begin{lemma}~\label{le33} 
Suppose that the BQDE system (\ref{eq201})--(\ref{eq203}) has an admissible solution
$\bx = (x_1 , x_2 )$.
Then one of the following holds.
\begin{description}
\item{(i)} The solution $\bx$ satisfies 
\beql{eq312}
 ||\bx|| < 200||F||^6 ~.
\eeq
\item{(ii)}
The equation
\beql{eq313}
y_1^2 - Dy_2^2 = g
\eeq
with $g$ given by (\ref{eq303}) has a solution $(y_1 , y_2 )$ such that
\beql{eq314}
y_1 > 0
\eeq
 and one of
 \beql{eq315a}
y_2 > 0~ \mbox{and ~$c(-b+ \sqrt D) > 0$}~, 
\eeq
\beql{eq315b}
y_2 < 0 ~\mbox{and ~$c(b+ \sqrt D) > 0 $}~, 
\eeq

 holds.
  In addition $(y_1 , y_2 )$ (mod $cD \Gamma )$ satisfies
\end{description}
\beql{eq316}
cy_1 + c(cd -be) \equiv cD \alpha_1 ( \bmod~ cD \Gamma ) ~,
\eeq
\beql{eq317}
-by_1 + Dy_2 +c(ae-bd) \equiv cD \alpha_2 ( \bmod ~cD \Gamma ) ~.
\eeq
 Conversely, if the system (\ref{eq313})-(\ref{eq317}) has a solution, then (\ref{eq201}) has
an admissible solution. 
\end{lemma}

\paragraph{Remark.} The lemma shows that either (i) the solution is small enough to
write down as a certificiate or (ii)  else we get definite  information on the sign conditions of the
solution values $(y_1, y_2)$ of the transformed equation (\ref{eq313}).

\paragraph{Proof.} 
Suppose that the admissible solution $\bx$ has
\beql{eq318}
|| \bx || \geq 200 ||F||^6 ~.
\eeq
We show that the solution $(y_1 , y_2 )$ to (\ref{eq313}) given by (\ref{eq301})
satisfies (\ref{eq314})--(3.17).
Now (\ref{eq305}) shows that the congruence conditions (\ref{eq316}), (3.17) hold.\\

We first show (\ref{eq318}) implies $|y_1 |$ and $|y_2 |$ are large.
If $|| \by || < 90 || F ||^5$ ,then absolute value estimates in (\ref{eq305}) yield
$$
|| \bx || < 190 || F ||^6
$$
contradicting (3.19).
So $|| \by || \geq 90 || F ||^5$.
If $|y_1 | < 90 || F ||^5$ then $|y_2 | \geq 90 || F ||^5$ and 
$$
y_1^2 \geq D(y_2)^2 - |g|
$$
implies
\beql{eq319}
|y_1 | > 89 || F ||^5 ~,
\eeq
so this holds in all cases.
Then
$$
y_2^2 \geq \df{y_1^2 - |g|}{D}
$$
implies
\beql{eq320}
|y_2| \geq \df{88}{\sqrt D} || F ||^5 ~.
\eeq

To prove (\ref{eq314}) holds, suppose for contradiction that $y_1 \le 0$.
Then by (\ref{eq301})
$$
y_1 - (be -cd) = Dx_1 \ge 0
$$
so
$$
|y_1 | \leq |be - cd | < 2 || F ||^2 ~,
$$
contradicting (\ref{eq319}).
Hence (\ref{eq314}) holds.\\

To prove one of (\ref{eq315a}) or (\ref{eq315b}) holds, note (\ref{eq313}) yields
\beql{eq321}
|y_1 + \sqrt D~ y_2  | ~ | y_1 - \sqrt D~ y_2  | = | g | \le  6 || F||^4 ~.
\eeq
Since
$$
y_1 + \sqrt D~ y_2  = (y_1 - \sqrt D~ y_2  ) + 2 \sqrt D~ y_2  ~,
$$
(\ref{eq320}) implies that
$$
MAX (|y_1 + \sqrt D~ y_2  |,~| y_1 - \sqrt D~ y_2  |) \geq 88 || F||^5 ~.
$$
Hence (3.21) yields
\beql{eq322}
MIN (|y_1 + \sqrt D~ y_2  |,~ |y_1 - \sqrt D~ y_2  |) \le  \frac{6}{88}||F||^{-1}< ||F||^{-1}.
\eeq
Consequently $y_1$ is very close to one of $\pm \sqrt D~ y_2 $.\\

We consider first the case that
\beql{eq323}
|y_1 - \sqrt D~ y_2  | < ||F||^{-1}~.
\eeq
Necessarily $y_2 > 0$, and
 (\ref{eq305}) yields
\begin{eqnarray}
cDx_2 & = & -by_1 + Dy_2 + c(ae - bd) \nonumber \\
 & = & \sqrt D (-b+\sqrt D) y_2 + \xi
\end{eqnarray}
where
\beql{eq325}
| \xi | \le |b| ||F||^{-1} + 2|| F ||^3 < 9 ||F||^4 ~.
\eeq
We claim
\beql{eq326}
|-b+c \sqrt D | \geq \df{1}{3} ||F||^{-1}~.
\eeq
Indeed, suppose $|-b+ \sqrt D | < 1$.
If so, then $|-b - \sqrt D | > 1$.
So
$$
|-b + \sqrt D | \geq \df{|b^2 -D|}{|-b- \sqrt D|} \geq \df{1}{3||F||}~,
$$
using $D < 2 ||F||^2$.
Now (\ref{eq320}) and (\ref{eq326}) imply
\beql{eq327}
| \sqrt D (-b+ \sqrt D) y_2 | > 29 ||F||^4 ~.
\eeq
Thus $\xi$ is too small to change
the sign of the two terms on opposite sides of (\ref{eq324}),
so that 
\beql{eq328}
\mbox{sign $(cDx_2) =$~sign $( \sqrt D (-b+\sqrt D )y_2 )$} ~.
\eeq
Since $x_2 \geq 0$, this yields
$$
\mbox{sign $(y_2 ) =$ ~sign $(c(-b+ \sqrt D ))$}
$$
which proves (\ref{eq315a}) holds in this case.
Next  we consider the case
$$
|y_1 + \sqrt D~ y_2 | < ||F||^{-1}~.
$$
Necessarily $y_2 < 0$ in this case, and 
an analysis similar to the previous case shows that
$$
\mbox{sign $(y_2) = - $~ sign $(c(b+ \sqrt D))$}
$$
which shows (\ref{eq315b})
holds in this case.\\

To prove the converse in (ii), suppose first 
that  we have a solution $(y_1 , y_2 )$ to conditions (\ref{eq313})--(\ref{eq317})
which satisfies (\ref{eq315a}).
Let $(t_1 , u_1 )$ be the minimal positive solution to Pell's equation
\beql{eq329}
t^2 - Du^2 = 1 ~.
\eeq
Let
\beql{eq330}
t_k + u_k \sqrt D = (t_1 + u_1 \sqrt D )^k ~.
\eeq
It is well-known that $(t_k , u_k )$ satisfy (\ref{eq329}) and that for any
modulus $M$ there exists an integer $P(M)$, the  period modulo M,  such that
\begin{eqnarray}~\label{eq331}
t_k & \equiv & 1 ~( \bmod ~M) \nonumber \\
u_k & \equiv & 0 ~( \bmod ~M)
\end{eqnarray}
whenever
$$
P(M)~ |~k ~.
$$
(See Appendix A.)
Certainly $t_k , u_k \rightarrow \infty$ as $k \rightarrow \infty$.
We now set
\beql{eq332}
y_1^\ast + y_2^\ast \sqrt D = (y_1 + y_2 \sqrt D) (t_k + u_k \sqrt D )
\eeq
where $P(cd \Gamma )~ | ~k$.
Note that since $y_1 , y_2 , t_k , u_k$ are all positive,
$y_1^\ast \geq t_k , y_2^\ast \geq u_k$.
By picking $k$ large enough, we may guarantee
\beql{eq333}
MIN (y_1^\ast , y_2^\ast ) > 88 ||F||^5 ~.
\eeq
Furthermore (3.31) applied with $M = cD \Gamma$ guarantees that
\beql{eq334}
y_i^\ast \equiv y_i ( \bmod ~cd \Gamma ) ~i = 1,2 ~.
\eeq
Also (\ref{eq332}) and (\ref{eq329}) guarantee that
$(y_1^\ast , y_2^\ast )$ satisfies the generalized Pell equation (\ref{eq313}).
Now let $(x_1^\ast , x_2^\ast )$ be the
rational solution to (\ref{eq201}) associated to
$(y_1^\ast , y_2^\ast )$ by (\ref{eq305}).
The congruence condition above implies that 
$(y_1^\ast , y_2^\ast )$ satisfies the congruences
(\ref{eq316}), (\ref{eq317}) hence
$(x_1^\ast , x_2^\ast )$ is an integer solution and
$$
x_i^\ast \equiv \alpha_i ( \bmod ~\Gamma ) ~.
$$
We claim that $(x_1^\ast , x_2^\ast )$ is nonnegative.
If so $(x_1^\ast , x_2^\ast )$ is the desired admissible solution.
Using (\ref{eq305}) and (\ref{eq333}), we obtain
$$
Dx_1^\ast
\geq y_1^\ast- |be - cd | > 86 ||F||^5,
$$
hence $x_1^\ast > 0$.
The bound (\ref{eq333}) implies that  the argument (\ref{eq321})--(\ref{eq322}) is valid,
 and since $y_1^\ast > 0$,
$y_2^\ast > 0$ we obtain
$$
|y_1^\ast - \sqrt D~ y_2^\ast | < ||F||^{-1} ~.
$$
The argument (3.24)--(3.28)  assumed only the truth
of (\ref{eq320}), so it is valid here as well and we obtain
$$
\mbox{sign $(cDx_2 ) =$ ~ sign $( \sqrt D (-b+c \sqrt D) y_2 )$}~.
$$
We are given $y_2 > 0$ by (\ref{eq315a}), and 
$c(-b+c \sqrt D) > 0$,
hence $x_2 > 0$ follows in this case.\\

Now assume a solution exists to (\ref{eq313})--(\ref{eq317}) which satisfies (\ref{eq315b}).
In this case set
$$
y_1^\ast + y_2^\ast \sqrt D = (y_1 + y_2 \sqrt D) (t_k - u_k \sqrt D)
$$
where $P(cD \Gamma ) ~|~k$.
Now  $y_1 > 0, y_2 < 0,$ and since $t_k > 0 , u_k > 0$, we obtain
$y_1^\ast > 0 , y_2^\ast < 0$ and
$y_1^\ast \geq t_k , |y_2^\ast | \geq |u_k |~$.
By picking $k$ large enough, we ensure that
\beql{eq335}
MIN (|y_1^\ast |,~|y_2^\ast | ) > 88 ||F||^5 ~.
\eeq
Also $y_i^\ast \equiv y_i$ (mod $cD \Gamma )$, and $(y_1^\ast , y_2^\ast )$
satisfies (\ref{eq313}).
An analogous argument to the case treated above now shows that
$(x_1^\ast , x_2^\ast )$ associated to this $(y_1^\ast , y_2^\ast )$ is an
admissible solution to (\ref{eq201}).
$~~~\Box$ \\

%
%
%
\subsection{Generalized Pell Equation: Exponential upper bound}

We next need an upper bound on
the size of admissible solutions to the generalized Pell equation
\beql{eq336}
y_1^2 - Dy_2^2 = g ~.
\eeq
The set of solutions to this equation has a simple form, related to solutions
of the Pell equation
$$
t^2 - Du^2 = 1 ~.
$$
Let $(t_1 , u_1 )$ denote the minimal positive solution to the Pell equation,
the
{\em fundamental solution}, and write 
\beql{eq337}
\epsilon = t_1 + u_1 \sqrt D ~.
\eeq
This is  a  unit in the real quadratic field $\QQ(\sqrt{D})$, which is 
the fundamental unit in the real quadratic order $\sO_D= \ZZ[1, \sqrt{D}]$. 
Note that $\bar{\epsilon}:= t_1 - u_1 \sqrt{D}= \epsilon^{-1}$ satisfies
$0 < \bar{\epsilon}  < 1.$
We recall the following upper bound on the size of the fundamental solution.

%
%

\begin{prop}~\label{pr35}
{\rm (Hua \cite{Hua42})}
Let $(t_1 , u_1 )$ be the minimal positive solution to Pell's equation
$t^2 - Du^2 = 1$.
If $\epsilon = t_1 + u_1 \sqrt D$ then
\beql{eq342}
\epsilon < D^{\sqrt D} ~.
\eeq
\end{prop}

Now we call a  solution $(y_1 , y_2 )$ to (\ref{eq336})  {\em basic}
provided $\eta = y_1 + y_2 \sqrt D$ has
\beql{eq338}
1 \leq | \eta | < \epsilon ~.
\eeq
the generalized Pell equation states  $\eta \bar{\eta}=g.$
We have the following finiteness result.

%
%
\begin{lemma}~\label{le34}
 For a positive squarefree $D$ the complete set of solutions  to
$$
y_1^2 - Dy_2^2 = g
$$
is given by $(y_1, y_2)= (y_{1,k}, y_{2,k})$ with
\beql{eq341}
y_{1,k} + y_{2,k} \sqrt D = \eta \epsilon^k
\eeq
for some basic solution $\eta$ and some integer $k$.
There are only a finite number of basic solutions.
\end{lemma} 

\paragraph{Proof.} 
Suppose $(y_1 , y_2 )$ is a solution to (\ref{eq336}).
Then for some integer $k$,
$$
\epsilon^k \leq |y_1 + y_2 \sqrt D | < \epsilon^{k+1}~.
$$
Consequently for the correct choice of sign
\beql{eq340}
\eta = \epsilon^{-k} (x_1 + x_2 \sqrt D )
\eeq
is a basic solution.

There are only a finite number of basic solutions since
$$
\bar{\eta} \equiv y_1 - y_2 \sqrt D = g/ \eta
$$
using (\ref{eq336}).
Hence
\begin{eqnarray}
|y_1 | & \leq & | \eta | + | \bar{\eta}| < \epsilon + |g| \nonumber \\
|y_2 | & \leq & \df{\epsilon + |g|}{\sqrt D} ,
\end{eqnarray}
as required.
$~~~\Box$ \\


We apply Lemma~\ref{le34} to establish 
an exponential  upper bound on the bit complexity  of 
writing down the least 
admissible solution to a genearalized
Pell equation,  if one exists,  in binary.

%
%
\begin{lemma}~\label{le36}
Suppose that the equation
\beql{eq343}
y_1^2 - Dy_2^2 = g
\eeq
has an integral solution $(y_1 , y_2 )$ with
$$
y_i \equiv \alpha_i ~( \bmod ~ M), ~ i = 1,2 ~,
$$
and with $(y_1 , y_2 )$ having prescribed signs.
Then it has such a solution with
\beql{eq344}
MAX ( \log |y_1| , \log |y_2| ) \leq 9 ||E||^{3/2} ( \log ||E||)^2
\eeq
where $||E|| = MAX (|D|, |g|, M)$. 
\end{lemma}

\paragraph{Proof.}
Consider the set of solutions $(y_{1,k} , y_{2,k} )$ to (\ref{eq343}) where
\beql{eq345}
y_{1,k} + y_{2,k} \sqrt D = \eta \epsilon^k ~,
\eeq
$k$ runs through the integers, and
$$
\eta = y_{1,0} + y_{2,0} \sqrt D
$$
where $(y_{1,0} , y_{2,0} )$ is a fixed basic solution of (\ref{eq343}).
We consider these solutions from the viewpoint of their sign patterns and
congruence classes (mod $M$).

For sign patterns, we will show that the signs of $(y_{1,k}, y_{2,k})$ 
become
constant for all sufficiently large positive $k$,
and also constant for  negative $k$ with $|k|$ sufficiently large. 
We show sign $(y_{1,k})$ is constant for all $k$ with
\beql{eq346}
k \geq \log |g| ~,
\eeq
and is constant for all $k$ with
\beql{eq347}
k \leq -\log |g| -2~.
\eeq
The same holds for $y_{2,k}$.
To do this, we use the standard notation $\bar{\alpha} = a - b \sqrt D$ 
for the
algebraic conjugate of $\alpha = a + b \sqrt D$.
Suppose (\ref{eq345}) holds.
Then
\beql{eq348}
y_{1,k} = \df{1}{2} ( \eta  \epsilon^k + \bar{\eta}  \bar{\epsilon}^k ) ~.
\eeq
Now $\bar{\eta} = g/ \eta$ by (\ref{eq343}) and
$\bar{\epsilon} = \epsilon^{-1}$.
Now suppose
\beql{eq349}
k \geq \log |g| \geq ( \log \epsilon )^{-1} \log |g|
\eeq
since the smallest $\epsilon$ that occurs is $\epsilon = 2 + \sqrt 3$
for $D = 3$, and so $\epsilon > e$.
Then since $1< \eta \le \epsilon$, $0< \bar{\epsilon}< 1$,
$$
| \eta \epsilon ^k | \geq |g| > | \df{g}{\eta}  \epsilon^{-1} |= 
| \bar{\eta} \bar{\epsilon}^k | ~.
$$
In this case $y_{1,k}$ has the same sign as $\eta$, and is constant.
Similarly when (\ref{eq347}) holds we find
$y_{1,k}$ has the same sign as $\bar{\eta}$.
(Use the fact $1 \leq | \eta | \leq \epsilon$.)
Similar arguments apply to $y_{2,k}$ using
$$
y_{2,k} = \df{1}{2 \sqrt D} ( \eta \epsilon^k - \bar{\eta}\bar{\epsilon}^k ) ~.
$$

For congruence conditions, in Appendix~A we show that for
\beql{eq350}
t_k + u_k \sqrt D = \epsilon^k
\eeq
and any modulus $M$, the sequences $\{t_k \}, \{u_k \}$ formed when $k$
varies are both periodic (mod $M$), and that the minimal positive period
$P(M)$ of both series jointly has
\beql{eq351}
P(M) \leq 2M( \log M+1 ) ~.
\eeq
Using (\ref{eq345}) this guarantees that the sequences
$\{y_{1,k} \}, \{y_{2,k} \}$ are both periodic (mod $M$) with
period $P(M)$.

Combining these results, we find that all possible combinations of sign
patterns and congruence conditions (mod $M$) that occur for
$(y_{1,k} , y_{2,k} )$ in the sequence (\ref{eq345}) occur for
some $k$ with
\beql{eq352}
|k| \leq 4M( \log M)+ \log g + 2 ~.
\eeq
In this circumstance
\begin{eqnarray*}
|y_{1,k} | & = & \df{1}{2} | \eta  \epsilon^k + \bar{\eta} \bar{\epsilon}^k | \\
~~~ \\
 & \leq & \df{1}{2} (\epsilon^{k+1} + |g| ) \leq |g| \epsilon^{k+1}
\end{eqnarray*}
using (\ref{eq338}),
$| \bar{\eta} | \leq |g|$, and $\bar{\epsilon} < 1$.
Hence
\begin{eqnarray}
\log |y_{1,k}| & \leq & \log g + (k+1 ) \log \epsilon \nonumber \\
 & \leq & (4M( \log M ) +2 \log g+3 ) \sqrt D \log D
\end{eqnarray}
using Proposition~\ref{pr35}.
Using
$$
|y_{2,k} | = \df{1}{2 \sqrt D} | \eta \epsilon^k - \bar{\eta}  \bar{\epsilon}^k |
$$
we obtain the same bound (3.53) for $\log y_{2,k}$.
The bound (3.53) implies
$$
MAX ( \log |y_1| , \log |y_2| ) \leq 9 || E ||^{3/2} ( \log ||E|| )^2 ~.
$$
By Lemma~\ref{le34} all integer solutions to (3.43) fall in one of the
sequences (\ref{eq345}), and this proves (\ref{eq344}).
$~~~\Box$ \\

%
%
%
\subsection{Indefinite BQDE:  Exponential upper bound}

Using Lemma \ref{le36}, 
we  obtain  an  upper bound on the bit complexity  
of an admissible solution $(x_1, x_2)$ to a general indefinite BQDE. 
In terms of the  number of binary digits
$\log |x_1| + \log |x_2|$ needed to write down this solution, it is singly exponential in terms of the
input size $L(F)$, since $L(F)$ is proportional to $\log ||F||$.
 However,   as a bound on the solution size $\max(|x_1|, |x_2|)$  it is   {\em double exponential} 
in terms of the input size $L(F)$.  

%
%
\begin{lemma}~\label{le37}
Any indefinite binary quadratic Diophantine equation that has an
admissible solution has such a solution $(x_1 , x_2 )$ with
\beql{eq354}
MAX ( \log |x_1| , \log |x_2| ) \leq 210 ||F||^6 ( \log ||F||)^2 ~.
\eeq
\end{lemma}

\paragraph{Proof.}
We apply the results of Lemma~\ref{le33} and Lemma~\ref{le36}.
Using (\ref{eq335}) a solution $(y_1^\ast , y_2^\ast )$ of
(\ref{eq343}) will correspond to an admissible solution
$(x_1^\ast , x_2^\ast )$ of (\ref{eq201}) provided (ii) of
Lemma~3.3 holds and
$$
MIN (|y_1^\ast |,~ |y_2^\ast |) \geq 88 || F ||^5 ~.
$$
By Lemma 3.4 we may write $(y_1^\ast , y_2^\ast ) = (y_{1,k} , y_{2,k} )$ for
some $\eta $ and $k$ in (\ref{eq345}).
But for
$$
k \geq 7 \log ||F|| + 11
$$
we have
\begin{eqnarray}
|y_{1,k} | & \geq & \df{1}{2} ( | \eta \epsilon^k | - | \bar{\eta} \bar{\epsilon}^k |) \nonumber \\
~~~ \\
 & \geq & \df{1}{2} [ \epsilon^k - |g|] \nonumber \\
~~~ \\
 & \geq &  e^{7 \log ||F||+10} - |g| \geq 88 ||F||^5 ~.
\end{eqnarray}
A similar bound holds for $|y_{2,k} |~$.
We obtain the same bound for
$$
k \leq -7 \log ||F|| - 13
$$
using
\begin{eqnarray*}
|y_{1,k} | & \geq & \df{1}{2} ( | \bar{\eta} \bar{\epsilon}^k | - | \eta \epsilon^k |) \\
~~~ \\
 & \geq & \df{1}{2} [  \epsilon^{|k|-1} -1 ] ~.
\end{eqnarray*}
and similarly for $y_{2,k}$.
Combining these inequalities with the argument of Lemma \ref{le36}, we find that if
(\ref{eq300}) has an admissible solution $(x_1 , x_2 )$ it has one whose
corresponding solution $(y_{1,k}, y_{2,k} )$ to (\ref{eq343}) (for some $\eta$) has
\begin{eqnarray*}
|k| & \leq & 4 |cD \Gamma |( \log cD \Gamma ) ~ 10 \log ||F|| + 20 \\
 & \leq & 90 ||F||^4 ( \log ||F||) ~.
\end{eqnarray*}
The same argument as Lemma \ref{le36} then gives
$$
MAX ( \log y_{1,k} , \log y_{2,k} ) \leq 100 ||F||^5
( \log ||F||)^2
$$
and (\ref{eq305}) then gives (\ref{eq354}).
$~~~\Box$
%
%
%

\section{Integral Binary Quadratic Forms}
\hsp
Section 3 essentially reduces the problem of finding succinct certificates
in the indefinite case to that of finding such
certificates for admissible  solutions of the generalized Pell equation
\beql{eq401}
y_1^2 - Dy_2^2 = g ~.
\eeq
The method of Gauss now relates solutions of this equation to 
the theory of integral binary forms; a solution to the equation will show
the equivalence of two particular binary forms, as given in section \S4.1.
Our problem will then be tranformed to finding succinct certificates establishing
such equivalence. \\

A final  step before applying the theory of integral binary forms is to reduce to the special 
case in which $y_1 , y_2$ are relatively prime; such $(y_1, y_2)$
are called {\em primitive solutions} to  (\ref{eq401}).
 If $y_1 , y_2$ is a solution to (\ref{eq401}) and
$(y_1 , y_2 ) = h$, set
\beql{eq402}
z_1 = \df{y_1}{h} , ~ z_2 = \df{y_2}{h} , ~
G = \df{g}{h^2}
\eeq
and obtain the equation
\begin{eqnarray}\label{eq403}
z_1^2  -  Dz_2^2 & = & G \nonumber \\
  (z_1 , z_2 ) & = & 1 ~.
\end{eqnarray}
If we specify  congruence conditions on $z_1 , z_2$ (mod $cD \Gamma )$ then we certainly
know $y_1 , y_2$ (mod $cD \Gamma )$.\\

An {\em integral binary quadratic form} $Q = [a,2b,c]$ is given by
$$
Q(x_1 , x_2 ) = ax_1^2 +2bx_1 x_2 + c x_2^2= \bx^T \bQ \bx~.
$$
Here
$$
\bx^T = [x_1 , x_2 ]
$$
and we call 
$$
\bQ = \left[
\begin{array}{c}
a~~ b \\
b~~c
\end{array}
\right]
$$
the symmetric matrix {\em associated to} the form $Q$.
The {\em determinant} $D$ of a form $Q = [a, 2b,c]$ is given by
\beql{eq405}
D = b^2 - ac = - \det(\bQ)~.
\eeq

A form $Q$ is {\em primitive} if the greatest common divisor $(a,b,c)=1$.
Primitive forms subdivide into {\em properly primitive forms} which have 
$(a ,2b,c) = 1$ and it  {\em improperly primitive forms} which have 
$(a, b, c)=1$ but $(a, 2b, c)= 2$.  
We shall  mainly deal with properly primitive forms in the
rest of this paper. \\

We say a form $Q$ {\em primitively represents} an integer $G$ provided
\beql{eq404}
Q(z_1 , z_2 ) = G
\eeq
for two relatively prime integers $z_1 , z_2$.
The {\em identity form} $I= I_D$ is $[1,0,-D]$, and it is properly primitive. 
In this terminology (\ref{eq403}) asserts that the identity form
primitively represents $G$.

%
%
%

\subsection{Equivalence of  Indefinite Binary Quadratic  Forms}

 Gauss  transformed the question of  (i) primitive representation of
an integer by a form to that of  (ii) determining the
equivalence of two forms.
A form $Q_1$ is {\em (properly) equivalent} to a form $Q_2$
if there is a $2 \times 2$ integer matrix $\bS \in SL(2, \ZZ)$ (i.e.  $\det ( \bS ) = 1$) such that
\beql{eq406}
\bS^T \bQ_1 \bS = \bQ_2 ~.
\eeq
This is an equivalence relation, and we denote it by $Q_1 \sim Q_2$, via $\bS$.
This equivalence relation preserves the determinant $D$, the property of
being a properly primitive form, and the property of primitively
representing a given integer $G$. 
%
%
\begin{lemma}~\label{le41}
 Let $z_1 , z_2$ be a primitive integer solution to  the generalized Pell equation
 $$
 z_1^2 - Dz_2^2=G.
 $$
 where $D$ is arbitrary.
 Then there exists  $z_3, z_4$ giving a  (proper) reduction matrix \\
  $\bZ = \left[
\begin{array}{c}
z_1 ~~ z_3 \\
z_2 ~~ z_4
\end{array}
\right]  \in SL(2, \ZZ)$
whose first column is $(z_1, z_2)^{T}$ 
which  shows the identity form \\
$I_D=[1,0, -D]$ of determinant $D$
is properly equivalent to a form  
\beql{eq407}
Q_0 = [G, 2B, C], 
\eeq
i.e $I \sim Q_0$, whose coefficients satisfy the bound
\beql{eq408}
MAX (2|B|, |C|) \leq |D| + 4G^2 ~.
\eeq
For any  choice of $B, C, z_3 , z_4$ satisfying (\ref{eq408}) we have
\begin{eqnarray}
MAX(|z_i |) & \leq & 6(|z_1 | + |D| +G^2 ) ~. \label{eq409}\\
MIN (|z_i |)  & \geq&  \df{1}{3 \sqrt D} (|z_1 | -5 |D| -5G^2 ) ~.\label{eq410}
\end{eqnarray}
\end{lemma}

\paragraph{Proof.}
Choose $z_3^\ast , z_4^\ast$ so that $z_1 z_4^\ast - z_2 z_3^\ast = 1$.
Then $\bS^\ast = \left[
\begin{array}{c}
z_1 ~~ z_3^\ast \\
z_2 ~~ z_4^\ast
\end{array}
\right]$
shows $I_D \sim Q^\ast$ where $I_D= [1, 0 . -D]$ and 
$$
Q^{\ast}= [G, 2B^\ast , C^\ast ],
$$
that is
$$
\bQ^{\ast} :=\left[\begin{array}{cc} 
G& B^{\ast}\\
B^{\ast}  & C^{\ast}\\
\end{array} \right]
 =( \bS^{\ast})^t \left[\begin{array}{cc} 
1 & 0 \\
0 & -D\\
\end{array} \right] \bS^{\ast}
$$
Now select $\lambda$ so that $\bS_1 := \left[
\begin{array}{c}
1~~\lambda \\
0 ~~ 1
\end{array}
\right]$ shows $Q^\ast \approx Q_0$ where
$$
Q_0 = [G, 2B,C] ,~ 0 < B < |G| ~.
$$
Then
$$
C = \df{D-B^2}{G}
$$
so
\beql{eq411}
|C| \leq |D| + G^2 ~,
\eeq
whence  (\ref{eq407}) and (\ref{eq408} ) hold. 
Also $I \sim Q_0$ via
$$
\bZ := \bS^\ast \bS_1 \equiv \left[
\begin{array}{c}
z_1 ~~ z_3 \\
z_2 ~~ z_4
\end{array}
\right] ~,
$$
where we have  $z_3= z_3^{\ast} + \lambda z_1, z_4= z_4^{\ast} + \lambda z_2$.\\

Now suppose that $\bZ$ is chosen to satisfy (4.7), (4.8) with  $I \sim Q_0$.
To bound the sizes of $z_2 , z_3 , z_4$ we observe that
\beql{eq412}
G = z_1^2 - Dz_2^2
\eeq
\beql{eq413}
B = z_1 z_3 - Dz_2 z_4
\eeq
\beql{eq414}
C = z_3^2 - Dz_4^2
\eeq
Then (\ref{eq412}) gives
\beql{eq415}
|| z_1 | - \sqrt D | z_2 || = \df{|G|}{|z_1 | + \sqrt D | z_2 |} \leq |G| ~.
\eeq
Hence
\beql{eq416}
\df{1}{\sqrt D} (|z_1 | - |G| ) < |z_2 | < \df{1}{\sqrt D} (|z_1 | + |G|)~.
\eeq
Similar arguments using (\ref{eq411}) and (\ref{eq414}) show
\beql{eq417}
||z_3 | - \sqrt D | z_4 || \leq \df{|C|}{|z_3 | + \sqrt D |z_4 |} \leq |C| \leq |D| + G^2,
\eeq
yielding
\beql{eq418}
\df{1}{\sqrt D} (|z_3 | - |D| -G^2 ) \leq |z_4 | \leq \df{1}{\sqrt D}
(|z_3 | + |D| +G^2 )~.
\eeq
Next note that
\begin{eqnarray}
B & = & \df{1}{2} \{(z_1 + \sqrt D z_2 ) (z_3 - \sqrt D z_4 ) + (z_1 - \sqrt D z_2 ) (z_3 + \sqrt D z_4 ) \} \nonumber \\
~~~ \nonumber \\
 & = & \df{1}{2} G \left\{ \df{z_1 + \sqrt D z_2}{z_3 + \sqrt D z_4} +
\df{z_3 + \sqrt D z_4}{z_1 + \sqrt D z_2 } \right\} ~.
\end{eqnarray}
Viewing this as $B = \frac{1}{2} G(x + \frac{1}{x} )$,  then $0 < B < |G|$
gives
$$
\df{1}{3} < |x| < 3 ~.
$$
If $z_3 , z_4$ have the same sign, then these bounds yield
\beql{eq420}
\df{1}{3} (|z_1 | + \sqrt D |z_2 |) < |z_3 | + \sqrt D |z_4 | < 3 (|z_1 | + \sqrt D |z_2 | ).
\eeq
Then using (\ref{eq416}) we obtain
\beql{eq421}
|z_3 | < 6 (|z_1 | + |G| )
\eeq
and
\beql{eq422}
|z_4 | < \df{6}{\sqrt D} (|z_1 | + |G|) ~.
\eeq
Combining (\ref{eq417}) and (\ref{eq420}), we obtain
$$
2|z_3 |  > \df{1}{3} (|z_1 | + \sqrt D |z_2 |) - (|D| +G^2)
$$
Now we can apply (\ref{eq415}) to obtain
 \nonumber \\
~\nonumber~~ \\
\beql{eqn423}
 2|z_3|  >  \df{1}{3} (|z_1 | - 3|D| - 3G^2 ) + \frac{1}{3}(|z_1|-|G| \ge  \df{1}{3} (2|z_1 | - 3|D| - 4G^2 ).
\eeq
Substituting this bound in the first inequality in  (\ref{eq418}) yields 
\beql{eq424}
|z_4 | > \df{1}{3 \sqrt D} (2|z_1 | - \frac{9}{2} |D| - 5 G^2 ) ~.
\eeq
If $z_3 , z_4$ have opposite signs, we use
$$
B = \df{1}{2} G \left\{ \df{z_1 - \sqrt D z_2}{z_3 - \sqrt D z_4} +
\df{z_3 - \sqrt D z_4}{z_1 - \sqrt D z_2} \right\}
$$
and again conclude the  bounds (\ref{eq420})--(\ref{eq424}) hold by similar arguments.
$~~~\Box$ \\

%
%
%

\subsection{Reduction of  Indefinite Binary Quadratic  Forms}

The problem of equivalence of indefinite forms, to determine if
$I_D \sim Q_0$ given by Lemma~\ref{le41}  may be simplified further.\\

Gauss introduced a notion of reduced indefinite form,
whose coefficients are bounded in absolute value b $2\sqrt{D}$, cf. (\ref{eq426}) below. 
He gave a reduction algorithm which shows that that each indefinite
form is properly equivalent to some reduced form.
This algorithm runs in polynomial time and is similar to the ordinary
continued fraction algorithm.
Application of this reduction algorithm
permits  the problem of equivalence of indefinite forms to be simplified to determining
equivalence of reduced indefinite forms.


%
%

\begin{defi}\label{de41}
{\em 
An  indefinite form $Q = [a,2b,c]$ is {\em reduced} when its coefficients
satisfy the bounds 
\begin{eqnarray}
0 & < & b < \sqrt D \nonumber\\
~~~~ \label{eq425}\\
\sqrt D - b & < & |a| < \sqrt D + b ~. \nonumber
\end{eqnarray}
}
\end{defi}

The reduction  inequalities  (\ref{eq425})  imply that any
(indefinite reduced form $Q_{red}$ satisfies 
$$
\sqrt D - b < |c| < \sqrt D + b
$$
so that
\beql{eq426}
| Q_{\rm red} | < 2 \sqrt D ~.
\eeq
There are in general many different reduced forms equivalent to a given
form; this is the subject of \S5. \\

Gauss's algorithm for reducing an indefinite form runs in polynomial
time, as given by the following bound.

%
%
\begin{prop}{\rm (Indefinite BQF Reduction Bound)}~\label{pr42}
 Given any indefinite form $Q$, there  exists a
reduced form $Q_{\rm red}$ and a reduction matrix $\bS_1 \in SL(2, \ZZ)$
such that $Q \sim Q_{\rm red}$ via $\bS_1$, and $\bS_1$ satisfies 
\beql{eq427}
\log || \bS_1 || = O ( \log || \bQ || ) ~.
\eeq
There is a reduction procedure which wehn given $Q$
will  obtain $Q_{red}$ and $\bS_1$ which
takes at most $O( \log ||\bQ|| M(\log(||\bQ||)$ elementary operations.
\end{prop}

\paragraph{Proof.} 
This bound is obtained in Lagarias \cite[Theorem 4.1]{Lag80a}.
 $~~~\Box$\\

We note that the (indefinite) identity form $I_D= [1, 0, -D]$ is not reduced in the sense above. 

\begin{defi} 
{\em  The {\em  reduced identity form}
$\tilde{I}:= \tilde{I}_D$ of positive, nonsquare  determinant $D >0$ is given by 
\beql{eq428}
\tilde{I}_D = [1,2 \lambda , \mu ]
\eeq
by $I _D \sim \tilde{I}_D$ via
\beql{eq429}
\bS^\ast = \left[
\begin{array}{c}
1 ~~ \lambda \\
0 ~~ 1
\end{array}
\right],~~~~ \mbox{with}~~~\lambda = \lfloor \sqrt{D} \rfloor,
\eeq
}
that is,  $(\bS^{*})^T \bQ_{I} \bS^{*}=\bQ_{\tilde{I}}$
and  $\mu=\lambda^2-D.$
\end{defi}

The reduced identity form $\tilde{I}_D$ is a reduced form in the sense above, and
is obtained by one step of the reduction algorithm in Proposition \ref{pr42}
applied to the identity form $I_D$. 

%
%
%

\subsection{Admissible Solutions of  Indefinite BQDE's and Reduced Forms}

The results obtained so far are summarized in the following lemma,
which will provide one part of the certificates. 
This lemma shows an equivalence between 
(i) existence of an admissible solution to
an indefinite BQDE (\ref{eq300}), and 
(ii)  equivalence of the  reduced identity form $\tilde{I}_D$ to a particular
reduced form $Q_{red}$ of determinant $D$, constructed using this BQDE.

%
%
\begin{lemma}~\label{le43}
 Consider the generalized Pell equation $E(y_1, y_2)=0$  given by
\beql{eq430}
 y_1^2 - Dy_2^2 = g ~,
\eeq
with $D>0$ not a perfect square.
This equation has a solution $(y_1 , y_2 )$ satisfying
\beql{eq431}
y_i \equiv \alpha_i ~( \bmod ~M) ,~~ i = 1,2,
\eeq
and with specified sign conditions
\beql{eq432}
sign (y_i ) = sign (i),~ i = 1,2,
\eeq
with $sign (1),~sign (2)$ given signs, if and only if  there exist
integers $h , B, C$ and $2 \times 2$ matrices $\bS , \bW \in SL(2, \ZZ)$ having the
following properties.
\begin{description}
\item{(i)}
 h is a positive integer and $G = g/h^2$ is an integer.
\item{(ii)}
The quadratic form $Q_0 = [G,2B,C]$ is properly primitive of determinant $D$.
\item{(iii)}
The matrix $\bS \in SL(2, \ZZ)$ shows
\beql{eq433}
Q_0 \sim Q_{\rm red} ~\mbox{via} ~ \bS
\eeq
where $Q_{\rm red}$ is a reduced form.
\item{(iv)}
 The matrix $\bW \in SL(2, \ZZ)$  shows
\beql{eq434}
\tilde{I}_D\sim Q_{\rm red} ~ \mbox{via} ~ \bW
\eeq
where $\tilde{I}_D$ is the reduced identity form and $Q_{\rm red}$ is given by (\ref{eq433}).
\item{(v)}
 Define $\bU = \left[ \begin{array}{c}
u_1 ~~ u_3 \\
u_2 ~~u_4
\end{array}\right] \in SL(2, \ZZ)
$
by
\beql{eq435}
\bU := \left[
\begin{array}{cr}
1 & -\lambda \\
0  & 1
\end{array}
\right]
\bW \bS^{-1}
\eeq
where $\lambda = \lceil \sqrt D \rceil$.
The congruence class of $\bW ( \bmod~M)$ is such that
\beql{eq436}
hu_i \equiv \alpha_i ( \bmod ~M), ~ i = 1,2 ~.
\eeq
In addition
\beql{eq437}
sign (u_i ) = sign (i), ~~ i = 1,2 ~.
\eeq
In fact $y_1 = hu_1 , y_2 = hu_2$ then satisfy
(\ref{eq430})--(\ref{eq432}).
If such an admissible solution exists, and we set
$$
||E||:= MAX (D, |g|, M)
$$
then there 
exist integers
$h, B, C$ and $2 \times 2$ matrices $\bS , \bW$ having properties
(i)--(v) and satisfying the bounds:
\item{(vi)}
\beql{eq438}
\log h = O ( \log ||E||) 
\eeq
\item{(vii)}
\beql{eq439}
MAX (|B|,|C|) \leq D + 4g^2 , 
\eeq
\item{(viii)}
\beql{eq440}
\log || \bS || = O ( \log ||E||), 
\eeq
\item{(ix)}
\beql{eq441}
\log || \bW || = O( ||E||^{3/2} ( \log ||E||)^2 )~. 
\eeq
\end{description}
\end{lemma}

\paragraph{Proof.}
Suppose that properties (i)--(v) above hold.
We check that $y_1 = hu_1 , y_2 = hu_2$ satisfy
(4.30)--(4.32).
The congruence and sign conditions hold by (4.36), (4.37), since $h > 0$
by (i).
To show (\ref{eq430}) holds, we observe that
$$
I_D \sim Q_0~ \mbox{via $\bU$}
$$
where $I_D= [1,0, -D]$ and $\bU$ is given by (\ref{eq438}).
For
$$
( \bS^\ast )^{-1} = \left[
\begin{array}{cr}
1 & - \lambda \\
0 & 1
\end{array}
\right]
$$
where $\bS^\ast$, given in (\ref{eq429}), shows $I_D \sim \tilde{I}_D$,
$\bW$ shows $\tilde{I}_D \sim Q_{\rm red}$ by (iv), and $\bS^{-1}$
shows $Q_{\rm red} \sim Q_0$ by (iii).
Thus
\beql{eq442}
\bu^t \left[
\begin{array}{cr}
1 & 0 \\
0 & -D
\end{array}
\right] \bU = \left[
\begin{array}{c}
G ~~ B \\
B ~~ C
\end{array}
\right] ~.
\eeq
Examining the upper left corner of this identity gives
$$
u_1^2 - Du_2^2 = G ~.
$$
Using $G = g/h^2 $ by (i), (\ref{eq430}) follows.

Now suppose an admissible solution $(y_1 , y_2 )$ satisfying (\ref{eq430})--(\ref{eq432}) exists.
Using Lemma~\ref{le36} we may suppose
\beql{eq443}
MAX ( \log |y_1| , \log |y_2| ) < 9 || E ||^{3/2} ( \log ||E||)^2 ~.
\eeq
Set
$$
h = g.c.d. (y_1 , y_2 ) ~.
$$
and $G= g/h^2$, establishing (i).
Since $h$ divides $g$,
$$
\log h = O( \log ||E||)
$$
giving the bound (vi).
Letting $z_1 = \frac{y_1}{h} , z_2 = \frac{y_2}{h}$, we may apply Lemma~4.1 to produce $B, C$
satisfying (ii) and the bound(vii), and a matrix $\bZ \in SL(2, \ZZ)$ showing
$I_D \sim Q_0$, and the lemma gives 
bounds (\ref{eq408}) yielding 
\beql{eq444}
\log || Q_0 || = O ( \log ||E||)
\eeq
and the bound (\ref{eq409}) gives
\beql{eq445}
\log || \bZ || = O (||E||^{3/2} ( \log ||E||)^2 )
\eeq
using (\ref{eq443}), since $z_1$ divides $y_1$.
Proposition~\ref{pr42} and (\ref{eq444}) produces an $\bS \in SL(2, \ZZ)$ satisfying (iii),
(viii).
Take $\bU = \bZ $ in (\ref{eq435}) and  use this equation to define $\bW$, namely
\beql{eq446}
\bW := \left[
\begin{array}{c}
1 ~~ \lambda \\
0 ~~ 1
\end{array}
\right] ~ \bZ \bS ~.
\eeq
The $u_1 = z_1 , u_2 = z_2$ so (v) holds.
Also $\left[ \begin{array}{c}
1 ~~ \lambda \\
0 ~~ 1
\end{array}
\right]$ shows $\tilde{I}_D \sim I_D, \bZ$ shows $I_D \sim Q_0$ and $\bS$
shows $Q_0 \approx Q_{\rm red}$,
hence (4.40) shows (iv) holds.
Finally (\ref{eq446}) gives
$$
|| \bW || \leq 8 || \bS^\ast || ~ || \bZ || ~ || \bS ||
$$
where $\bS^\ast = \left[ \begin{array}{c}
1 ~~ \lambda \\
0 ~~ 1
\end{array}
\right]$ with $\lambda = \lceil \sqrt D \rceil$.
Hence
$$
\log || \bW || = O ( ||E||^{3/2} ( \log ||E||)^2 )
$$
using (\ref{eq445}) and the already established (viii).
$~~~\Box$
%
%
%

\section{Equivalence of Reduced Indefinite Binary Quadratic Forms}
\hsp
The problem has now been simplified to that of finding a particular matrix $\bW$
which demonstrates the equivalence of the reduced identity form $\tilde{I}_D$ and a
reduced form $Q_{\rm red}$.
Gauss [11, Arts. 183--205] showed that the set of all reduced indefinite
forms  of determinant $D$ has a simple
structure, which we describe below.

%
%
%
\subsection{Cycles of Indefinite Reduced Forms and the Principal Cycle}

An (indefinite) reduced form $Q_1 = [a_1 , 2b_1 , c_1 ]$ is said to have as a {\em right neighbor}
the reduced form 
$Q_2 = [a_2 , 2b_2 , c_2 ]$  provided $a_2 = c_1$.
In this case $Q_2$ is unique and $Q_1 \sim Q_2$ via $\bS$ where
\beql{eq501}
\bS = \left[
\begin{array}{cl}
0 & 1 \\
-1 & \lambda
\end{array}
\right],
\eeq
in which  $\lambda$ is specified by
\beql{eq502}
- \sqrt D - b_1 < \lambda c_1 < - \sqrt D - b_1 + | c_1 | ~.
\eeq
Travelling to right neighbors results in traversing a cycle of reduced forms.
The collection of all reduced indefinite primitive forms of determinant $D$ (which is finite)
partitions  into a finite set of cycles of possibly different lengths under the right-neighbor
relation.\\

The cycle containing the reduced principal form $\tilde{I}_D$ is called the {\em principal cycle}.
Let $Q^{(1)} $ denote the right-neighbor of $\tilde{I}_D$, and
$Q^{(j)}$ the right-neighbor of $Q^{(j-1)}$.
Let $\bS^{(j)}$ denote the matrix given by (\ref{eq501}), (\ref{eq502}) taking
$Q^{(j-1)}$ to $Q^{(j)}$.
The set of $Q^{(j)}$ form a closed {\em cycle} of even {\em period}
$2p$, i.e., there exists some $Q^{(k)} = \tilde{I}$ and the
smallest $k = 2p$.
For $1 \leq j \leq 2p$, $\tilde{I} \sim Q^{(j)}$ via $\bL_j$ where
\beql{eq503}
\bL_j = \bS^{(1)} \ldots \bS^{(j)} ~.
\eeq
We call $\bL_j$ a {\em simple equivalence matrix.}
The matrix $\bU = \bL_{2p}$ is called the
{\em fundamental automorph}.
If we set
$$
\bU = \left[
\begin{array}{c}
u ~~ w \\
t ~~ v
\end{array}
\right]
$$
the condition $\tilde{I} \sim \tilde{I}$ via $\bU$ shows that $t,u$
satisfies Pell's equation
\beql{eq504}
t^2 - Du^2 = 1 ~.
\eeq
In fact $(|t|,|u|)$ is the least strictly positive solution to (\ref{eq504}),
the
{\em fundamental solution}, and
\beql{eq505}
\bU = \left[
\begin{array}{ll}
u & -t  \\
t & -Du
\end{array}
\right] ~.
\eeq
We may consistently extend the definition of $\bL_j$ to apply for all
integers $j$ by first defining $\bS^{(j)}$ for negative $j$ by
$$
\bS^{(j)} = [ \bS^{(j_0 )} ]^{-1}
$$
where $j \equiv j_0 ( \bmod ~ 2p )$ and $0 < j_0 \leq 2p$, using (\ref{eq503})
for all positive $j$ and using
$$
\bL_{-j} = \bS^{(-j)} \ldots \bS^{(-1)}
$$
for $j > 0$.
In that case, for any integer $k$ we have
$$
\bL_{j+2kp} = \bU^k \bL_j ~.
$$
Gauss proved the following result (see Mathews \cite[Arts. 76, 88]{Mat61}, Venkov \cite{Ven70}). \\
%
%
\begin{prop}~\label{pr51}
{\rm (Gauss)} Let $\tilde{I} \sim Q$ via $\bT$ where $Q$ is reduced.
Then there is some $j$ with $1 \leq j \leq 2p$
such that $Q = Q^{(j)}$.
Furthermore there is an integer $k$ such that
\beql{eq507}
\bT = \pm \bU^k \bL_j = \pm \bL_{j+2kp} ~.
\eeq
\end{prop}

%
%
%
\subsection{Sign Patterns of Equivalence Matrices $\bL_j$ }

To handle the nonnegativity conditions
in  Theorem~\ref{th11},  we will need  detailed information about the
signs of entries in the equivalence matrices $\bL_j$.
We first introduce the notation that if a matrix $\bM = [m_{ij}]$, then
\beql{eq508}
| \bM | = [|m_{ij}|]~.
\eeq

%
%
\begin{lemma}~\label{le52}
 The equivalence matrices $\bL_j$ have the following properties.
\begin{description}
\item{(i)}
 For $j > 0$ the entries of $L_j$ have the sign patterns
$\left[ \begin{array}{c}
++ \\
++
\end{array}
\right]$,
$\left[ \begin{array}{c}
-+ \\
-+
\end{array}
\right]$,
$\left[ \begin{array}{c}
-- \\
--
\end{array}
\right]$,
$\left[
\begin{array}{c}
+- \\
+-
\end{array}
\right]$
 according as $j = 0, 1, 2$ or $3 ~(\bmod~ 4)$.
\item{(ii)}
{\em For $j > 0$,}
\beql{eq509}
| \bL_j | = | \bS^{(1)} | \ldots | \bS^{(j)} |
\eeq
 and
\beql{eq510}
| \bL_{-j} | = | \bS^{(-j)} | \ldots | \bS^{- (1)} | ~.
\eeq
\item{(iii)}
 The four entries of $\bL_j = ( l_{ij})$ are all about the same size in the
sense that for any $|j| \geq 2$,
\beql{eq511}
1 \leq \df{MAX |l_{ij} |}{MIN | l_{ij}|} \leq 4(D+ \sqrt D) ~.
\eeq
\end{description}
\end{lemma}

\paragraph{Proof.} 
We first observe that a reduced form $Q = [a,2b,c]$ has by definition
(\ref{eq420}) $|b| < \sqrt D$ hence
\beql{eq512}
ac < 0 ~.
\eeq
The reduced forms $Q^{(i)} = [a_i , 2b_i , c_i ]$ in the fundamental cycle
have $a_{i+1} = c_i$.
Noting $Q^{(0)} = \tilde{I}$ so
$a_0 = 1$, by induction using (\ref{eq512}) we obtain
\beql{eq513}
(-1)^i a_i > 0 ~.
\eeq
Now (4.25) and (\ref{eq502}) imply that
$$
\lambda_i c_{i-1} < 0 ~,
$$
so that (\ref{eq512}), (\ref{eq513}) yield
\beql{eq514}
(-1)^{i+1} \lambda_i > 0 ~.
\eeq

To prove (i) and (ii), note (\ref{eq514}) implies $\bS^{(i)}$ for
$i > 0$ has the sign patterns $\left[ \begin{array}{c}
-+ \\
-+
\end{array}
\right]$ when $i$ is odd,
$\left[ \begin{array}{c}
++ \\
--
\end{array}
\right]$ when $i$ is even.
Then it is easy to establish by induction on $i > 0$ that the entries of
$\bL_i$, have the sign patterns
$\left[ \begin{array}{c}
++ \\
++
\end{array}
\right]$,
$\left[ \begin{array}{c}
-+ \\
-+
\end{array}
\right]$,
$\left[ \begin{array}{c}
-- \\
--
\end{array}
\right]$,
$\left[ \begin{array}{c}
+- \\
+-
\end{array} \right]$
according as $i = 0, 1, 2, 3$ (mod 4).
Another induction on $i > 0$ shows no cancellation occurs in
multiplying the entries of $\bL_i$ by $\bS^{(i+1)}$ and (\ref{eq509}) follows.
For $i < 0$, we observe first that
\begin{eqnarray*}
\left[
\begin{array}{cl}
0 & l \\
-1 & \lambda
\end{array}
\right]^{-1} = \left[
\begin{array}{lr}
\lambda & -1 \\
1 & 0
\end{array}
\right] ~.
\end{eqnarray*}
Then note $\lambda_i = \lambda_{i-2p}$ so (\ref{eq514}) holds for $i < 0$ as well.
This implies that for $i < 0$, $\bS^{(i)}$ has the sign patterns
$\left[ \begin{array}{c}
-- \\
++
\end{array}
\right]$ if $i$ is odd,
$\left[ \begin{array}{c}
+- \\
+-
\end{array}
\right]$ if $i$ is odd,
$\left[ \begin{array}{c}
+- \\
+-
\end{array}
\right]$ if $i$ is even.
Another induction shows for $i < 0$ that the entries of $\bL_i$ have the
sign patterns
$\left[ \begin{array}{c}
++ \\
++
\end{array}
\right]$,
$\left[ \begin{array}{c}
-- \\
++
\end{array}
\right]$,
$\left[ \begin{array}{c}
-- \\
--
\end{array}
\right]$,
$\left[ \begin{array}{c}
++ \\
--
\end{array}
\right]$ according as$i \equiv 0,1,2$ or 3 (mod 4).
Then (\ref{eq510}) follows by induction.

To prove (iii), consider first the case $j > 0$.
Using (\ref{eq509}), we need only bound the entries
\beql{eq515}
\bL_j | = \left| \begin{array}{lc}
0 & 1 \\
1 & | \lambda_1 |
\end{array}
\right| \ldots \left|
\begin{array}{lc}
0 & 1 \\
1 & |\lambda_j |
\end{array}
\right| ~.
\eeq
The formulae for $| \bL_j |$ is exactly that of the ordinary continued
fraction algorithm, where
\beql{eq516}
| \bL_j | = \left| \begin{array}{c}
p_{j-1}~~~ p_j \\
q_{j-1} ~~~ q_j
\end{array} \right|
\eeq
and $\frac{p_j}{q_j}$ is the $j$th convergent to $\theta = [0, | \lambda_1 |, | \lambda_2 | , \ldots ]$.
In particular, for any $j \geq 2$ we have
\beql{eq517}
\df{| \lambda_2 |}{|\lambda_1 \lambda_2 | +1} = \df{p_3}{q_3} \leq \df{p_j}{q_j} \leq \df{p_2}{q_2} = \df{1}{| \lambda_1 |} ~.
\eeq
Now (\ref{eq502}) implies
\beql{eq518}
1 \leq | \lambda_i | < 2 \sqrt D
\eeq
so (\ref{eq516}) yields
\beql{eq519}
\df{1}{2 \sqrt D + 1} \leq \df{p_j}{q_j} < 1 ~.
\eeq
In addition
\beql{eq520}
q_{j+1} = | \lambda_j | q_j + q_{j-1} \leq (| \lambda_j | +1) q_j~.
\eeq
Combining (\ref{eq517}), (\ref{eq518}), we obtain
\begin{eqnarray}
p_j & \leq & p_{j+1} \leq q_{j+1} \nonumber \\
p_j & \leq & q_j \leq q_{j+1} ~.
\end{eqnarray}
Finally
\beql{eq522}
q_{j+1} \leq (2 \sqrt D + 1) q_j \leq (2 \sqrt D + 1)^2 p_j
\eeq
and since $p_j \geq 1$ for $j \geq 2$ by (\ref{eq516}) this implies (\ref{eq511})
on this range.
The case $j < 0$ is treated analogously to $j > 0$.
In this case
$$
| \bL_{-j} | = \left|
\begin{array}{c}
q_{j+1} ~~ q_j \\
p_{j+1} ~~ p_j
\end{array}
\right| ~,
$$
however.
$~~~\Box$ \\

\paragraph{Remark.}  There is a close connection between the $\lambda_i$ and the
ordinary continued fraction (OCF) expansion of $\sqrt D$.
It is known that the OCF expansion has the form 
$$
\sqrt D = [ \mu_0 , \overline{\mu_1 , \ldots , \mu_n} ]
$$
in which  $[\mu_1 , \ldots , \mu_n ]$ is the purely periodic part of the expansion
and $n$ is the shortest period. (Stark \cite[Sec. 7.7]{St70}). The $2 \times 2$ matrices in the
continued fraction expansion have determinant $-1$ (see Stark    \cite[Sec. 7.6]{St70}) 
while the matrices in visiting neighboring forms in the 
Gaussian reduction procedure have determinant $+1$, but the entries of 
the resulting two expansions  are simply related up to signs, and one obtains that 
$n = 2p$ if $n$ is even, and $n = p$ otherwise.
In either case $\mu_i = | \lambda_i |$ for $1 \leq i \leq 2p$.

%
%
%
\subsection{Bounds on  Sizes of Equivalence Matrices $||\bL_j||$}

Our next step is to  estimate the size of the entries of $\bL_j$ in relation to $j$. 

%
%
\begin{lemma}~\label{le53}
 For all $j > 0$,
\beql{eq523}
\log || \bL_{j+2} || \geq \log || \bL_j || + 1 ~,
\eeq
 and
\beql{eq524}
\log || \bL_j || \leq \log || \bL_{j+1} || \leq \log || \bL_j || + \log D+2 ~.
\eeq
\end{lemma}

\paragraph{Proof.}
For the bound (\ref{eq523}), we use (\ref{eq515}), (\ref{eq516}) to obtain
$$
q_{j+2} = (| \lambda_{j+1} \lambda_j | + 1) q_j + | \lambda_{j+1} | q_{j-1} \geq 2q_j ~.
$$
The left side of (\ref{eq524}) follows from (5.21).
Finally (\ref{eq522}) implies
$$
\log || \bL_{j+1} || \leq \log || \bL_j || + \log (2 \sqrt D + 1 )
$$
from which the right side of (\ref{eq524}) follows.

Analogous inequalities to (\ref{eq523}), (\ref{eq524}) hold for $j < 0$.$~~~\Box$ \\

Using Lemma~\ref{le53} it is easy to prove, by induction that for $j \geq 1$, that we have
\beql{eq525}
\df{1}{2} |j| \leq \log || \bL_j || \leq |j| ( \log D+2 ) ~.
\eeq
The same holds for $j \leq -1$.\\

%
%
%
\subsection{Recursion Formula for $\bW$ in Terms of $\bL_j$'s}\label{sec54}

We can use the preceding results to determine a matrix formula 
for the equivalence matrix $\bW$ of Lemma~\ref{le43}. \\
%
%
\begin{lemma}~\label{le55}
Suppose that the indefinite BQDE $y_1^2-Dy_2^2= g$
in Lemma~\ref{le43}  has an admissible solution.
 Let $\bW$ be the equivalence matrix guaranteed to exist in Lemma~\ref{le43},
satisfying (i)--(ix) of that lemma, and certifying an admissible solution.
Then
\beql{eq528}
\bW = \pm ( \bL_{2p})^k \bL_j
\eeq
 for some $j$ with $1 \leq j \leq 2p$ and an integer $k$ satisfying
\beql{eq529}
|k| = O ( ||E||^{3/2} ( \log ||E||)^2 ),
\eeq
with $||E||= MAX(D, |g|, M).$
\end{lemma} 

{\bf Proof.}
This follows immediately on combining Proposition~\ref{pr51}, Lemma~\ref{le53} 
with  the size bound (4.41) on $\bW$. The role of the extra power $k$ is to 
meet the side congruence conditions.
$~~~\Box$

%
%
%
\subsection{Upper Bound on Length of the Principal Cycle.}
We next give an upper bound for the length $2p$ of the principal cycle. 
This upper bound holds more generally to all cycles of reduced forms of determinant $D$.   \\
%
%
\begin{prop}~\label{pr54}
 The period $2p$ of the fundamental cycle of reduced forms of 
 positive nonsquare determinant $D$ satisfies
\beql{eq526}
p < ( \sqrt D + 1) \log D.
\eeq
\end{prop}
\paragraph{Proof.}
The result of Hua \cite{Hua42} given in Proposition \ref{pr35}
asserts that if $(t_0 , u_0 )$ is the fundamental
positive solution to $x^2 - Dy^2 = 1$ then
$$
\df{t_0 + u_0 \sqrt D}{2} < D^{\sqrt D} ~.
$$
Using (5.5) we obtain
\beql{eq527}
|| \bL_{2p} || \leq D^{\sqrt D +1} ~.
\eeq
Combining (\ref{eq527}) with (\ref{eq525}) gives
$$
p \leq \log || \bL_{2p} || \leq ( \sqrt D + 1) \log D ~. ~~~\Box
$$

\paragraph{Remark.} The examples $D = 5^{2n+1}$ with period $p=5^n$
mentioned in the introduction show that periods
$p > \frac{1}{3} \sqrt D$ do occur.

%
%
%
\subsection{Exponential Time Algorithm for Solving a BQDE}
\label{sec56}

We now consider the general binary quadratic Diophantine equation (\ref{eq300})
in standard form, 
$$
ax_1^2 + 2 b x_1x_2 + cx_2^2 + 2dx_1 + 2ex_2 + f=0.
$$
but with no side conditions imposed. We give an exponential time running
bound for determining if the equation has an integer solution. The algorithm
analyzed is a variant of the method of Gauss to find an integer solution. 
The algorithm is simplified since our 
object is only to obtain a bound of form $O\left( 2^{c_1 L(F)}\right)$,without
optimizing the constant $c_1$.

%
%
\paragraph{Proof of Theorem~\ref{th14}.}
We reduce the Diophantine equation in Theorem~\ref{th14} to the standard
form by multiplying its coefficients by $2$ if necessary. We consider the
following algorithm.  If the BQDE is definite or degenerate, it suffices
by Lemma~\ref{le32} to sequentially test all integer vectors $(x_1, x_2)$ with $||\bx|| \le 8||F||^4$
to see if they satisfy the equation. This takes at most 
$O\left( ||F||^4 M( \log ||F||)\right) = O\left( 2^{c_0^{\ast} L(F)}\right)$
elementary
operations. If the BQDE is indefinite, we reduce it in polynomial time to a generalized Pell equation
$E(y_1, y_2)=0$ with 
$E(y_1, y_2) = y_1^2- D y_2^2-g$, as in Lemma~\ref{le32},
noting that $\log ||E|| = O\left(L(F)\right).$ \\

We solve the indefinite case by checking all possible certificates for solutions that are of
the form given by Lemma~\ref{le43} (i)-(ix).
To do this  we first find all square divisors $h^2$ of $g$ by exhaustive
search, set $G= \frac{g}{h^2}$ and determine
all properly primitive candidate forms $Q_0=[G, B, C]$ of
determinant $D$ having $MAX(|B|, |C|) \le D+ 4g^2$. This can be done by
enumeration in 
$O\left( (D+4g^2) g\right) =O\left(2^{c_1^{\ast}L(F)}\right)$ 
elementary operations. For each such quadratic form $Q_0$, an integer solution to the
BQDE will exist if it is equivalent to $\tilde{I}_D$. If none of the forms $Q_0$ are equivalent to the 
reduced identity form 
$\tilde{I}_D$, then the BQDE has no integer solution.\\

To test if a given indefinite form $Q_0$ is  equivalent  to $\tilde{I}_D$, we first reduce 
$Q_0$ to an indefinite reduced form $Q_{red}$, which by Proposition \ref{pr42} takes 
at most $O\left( \log ||Q_0||) M(\log||Q_0||)\right)= O \left( \log ||E||)^3 \right)$ bit  operations.
Next one tests if $Q_{red}$  is one
of the $2p$ forms in the principal cycle. To do this it suffices to step through all
forms in the principal cycle, starting with $\tilde{I}_d$ and see if there is a match.
By Proposition~\ref{pr54} there are at most $\sqrt{D+1} \log D$ forms in the
cycle, and this test takes at most $O\left(\sqrt{D}(\log D)^3\right)$ bit operations.
We conclude that  a single such test takes at most $O\left(2^{c_2^{\ast} L(F)}\right)$ 
elementary  operations.
 Combining all these tests, we can (wastefully) take 
 the constant $c_1= c_0^{\ast} +c_1^{\ast}+ c_2^{\ast}$ to get a
running time bound $O\left( 2^{c_1 L(F)} \right)$ elementary operations on the algorithm.
$~~~\Box$

%
%
%

\section{Composition of Binary Quadratic Forms and Infrastructure}
\hsp
The certificates we construct use the operation of composition of binary
quadratic forms introduced by Gauss \cite{Gau66}, in particular the action of composition
on the fundamental cycle of reduced forms.
Our treatment of composition of forms is based on Mathews \cite{Mat61}, Venkov \cite{Ven70}
and Lagarias \cite{Lag80a}.
One may also consult   Buell \cite{Bu89}, 
Shanks \cite{Sha71}, and Smith \cite[Arts. 105--113]{Smi}.\\
 
 The idea of analyzing the action of composition on the cycle of reduced
forms equivalent to $\tilde{I}$   is due to D.~Shanks \cite{Sha72},
whose called it the ``infrastructure".
We give bounds on the infrastructure in terms of composition of forms.

%
%
%

\subsection{Composition of Binary Quadratic Forms}

The simplest example of composition of two binary quadratic forms is the identity
\beql{eq601}
(x_1^2 + x_2^2 ) (y_1^2 + y_2^2 ) = (x_1 y_1 -x_2 y_2 )^2 + (x_1 y_2 + x_2 y_1 )^2
\eeq
noted by Fermat. This identity  shows that the product of two numbers which are the
sum of two squares is itself the sum of two squares.
We can rewrite (\ref{eq601}) in the form
\beql{eq602}
Q_1 ( x_1 , x_2 ) Q_1 (y_1 , y_2 ) = Q_1 (x_1 y_1 - x_2 y_2 , x_1 y_2 + x_2 y_1 )
\eeq
where $Q_1 = [1,0,1]$, and in matrix terms as
\beql{eq603}
(\bx^T \bQ_1 \bx)( \by^T \bQ_1 \by) = \bz^T\bB^T \bQ_1 \bB {\bz}
\eeq
where
\begin{eqnarray}\label{eq604}
\bx^T& = & [x_1 , x_2 ] ,~~~~~\by^T = [y_1 , y_2 ] \nonumber \\
\bz^T& = & [x_1 y_1 , x_1 y_2 , x_2 y_1 , x_2 y_2 ]
\end{eqnarray}
and
\beql{eq605}
\bB = \left[
\begin{array}{lllr}
1 & 0 & 0 & -1 \\
0 & 1 & 1 & 0
\end{array}
\right] ~.
\eeq
In this case we say {\em  $Q_1$ is composed of $Q_1$ and $Q_1$ via the bilinear matrix
$\bB$} of (\ref{eq605}).\\

In the general case we say a quadratic form $Q_3 = [a_3 , 2b_c , c_3 ]$ is
{\em composed} of forms $Q_1 = [a_1 , 2b_1 , c_1 ]$ and $Q_2 = [a_2 , 2b_2 , c_2 ]$ via a {\em bilinear matrix} $\bB$ provided the matrix equation
\beql{eq606}
\bx^t \bQ_1 \bx \by^t \bQ_2 \by = \bz^t \bB^t \bQ_3 \bB \bz
\eeq
holds, where $\bx , \by , \bz$ are given by (\ref{eq604}), the $x_i$ and $y_j$
are indeterminates.
Here $\bB$ is a $2 \times 4$ integer matrix (``bilinear matrix"), which is  required to be
{\em unimodular} and {\em oriented}
(terms defined below). We write $Q_3 = Q_1 \circ Q_2$ to indicate composition of forms,
with the associated bilinear matrix $\bB$ being omitted from the notation. \\

%
%

\begin{defi}\label{de61}
{\em 
A $2 \times 4$ integer matrix $\bB = [b_{ij} ]$ is said to be:
\begin{description}
\item{(i)}
{\em unimodular} provided the six cofactors
$$
\Delta_{ij} = \left[ \begin{array}{c}
b_{1i} ~~~ b_{1i} \\
b_{2i} ~~~ b_{2j}
\end{array}
\right] ,~~ 1 \leq i < j \leq 4
$$
have greatest common divisor 1.
\item{(ii)}
{\em oriented}
provided $a_1 \Delta_{12} > 0$ and $a_2 \Delta_{13} > 0$.
\end{description}
}
\end{defi}

%
%

Recall that a form $Q=[a,2b,c]$ is {\em properly primitive} if  $GCD(a, 2b, c)=1$.
If $Q_3$ is the composition of two properly primitive forms $Q_1$ and $Q_2$ of
determinant $D$, then $Q_3$ itself is properly primitive of determinant $D$.
This is a consequence of the unimodularity property of $\bB$.\\

In the rest of this section we deal only with {\em properly primitive}
indefinite forms. We use the following result on  composition.
%
%
\begin{prop}~\label{pr61}
 Given any two properly primitive  reduced forms $Q_1 , Q_2$ of determinant $D$ there is a
properly primitive reduced form $Q_3$ of determinant $D$ and a bilinear matrix $\bB$ such that
$Q_3 = Q_1 \circ Q_2$ via $\bB$ and
\beql{eq607}
\log || \bB || = O( \log D ) ~.
\eeq
There is an algorithm which when given as input $Q_1, Q_2$ in binary will
determine  $Q_3$ and matrix $ \bB$ in binary, which runs in 
at most $O \left( (\log D) M(\log D)\right)$ bit operations.
\end{prop}

\paragraph{Proof.} This is shown in Lagarias \cite[Theorem 5.5]{Lag80a}. $~~~\Box$\\

%
%
%

\subsection{Infrastructure Bounds}

The key result facilitating the use of composition to create short certificates
is the following lemma.
Before stating it, we recall that the Kronecker product $\bS \otimes \bT$
of an $m \times n$ matrix $\bS = [s_{ij} ]$ and a $k \times l$ matrix $\bG$
is a $km \times ln$ matrix
$$
\bS \otimes \bT = \left[
\begin{array}{c}
s_{11} \bT \\
\vdots \\
s_{m1} \bT
\end{array}
\begin{array}{c}
\ldots \\
~~~ \\
\ldots
\end{array}
\begin{array}{c}
s_{1n} \bT \\
\vdots \\
s_{mn} \bT
\end{array}
\right]
$$
given in block matrix form.

%
%
\begin{lemma} ~\label{le62}
 Let $\tilde{I}_D \sim Q_1$ via $\bS_1$ and $\tilde{I}_D \sim Q_2$ via $\bS_2$.
If $Q_3 = Q_1 \circ Q_2$ via $\bB$, then $\tilde{I}_D \sim Q_3$ via $\bS_3$
where $\bS_3$ satisfies the matrix equation
\beql{eq608}
\bS_3 \bB = \bB_0 ( \bS_1 \otimes \bS_2 )
\eeq
where
\beql{eq609}
\bB_0 = \left[
\begin{array}{c}
1~~0~~0 \\
0~~1~~1
\end{array}
\begin{array}{c}
D- \lambda^2 \\
2 \lambda
\end{array}
\right]
\eeq
and $\lambda = [ \sqrt D ]$. 
\end{lemma}

\paragraph{Proof.}
It is straightforward to check that the identity form $I_D = I_D\circ I_D$ via
$$
\bB = \left[
\begin{array}{c}
1 ~~ 0 ~~ 0 ~~D \\
0 ~~ 1 ~~ 1 ~~ 0
\end{array}
\right] ~.
$$
Using Lagarias \cite[Lemma 5.1 (i)]{Lag80a}, since $I_D \sim \tilde{I}_D$
via $\left[ \begin{array}{c}
1 ~~ \lambda \\
0 ~~ 1
\end{array}
\right]$ 
we obtain that the reduced identity form has $\tilde{I}_D = \tilde{I}_D \circ \tilde{I}_D$ via $\bB_0$.
Using the same \cite[Lemma 5.1 (i)]{Lag80a},  we next conclude $\tilde{I}_D = Q_1 \circ Q_2$
via $\bB_0 ( \bS_1 \oplus \bS_2 )$.
Then using \cite[Lemma 5.1 (ii)]{Lag80a}, we conclude there exists an
integer matrix $\bS_3 \in SL(2, \ZZ) $ such that
$$
\bS_3 \bB = \bB_0 ( \bS_1 \otimes \bS_2 ) ~,
$$
the desired result.
$~~~\Box$\\

We note that $\bS_3$ is uniquely determined by equation (\ref{eq608}), since $\bB$
contains an invertible $2 \times 2$ submatrix by the unimodularity condition.\\

Now suppose $Q_1$ and $Q_2$ are forms in the principal cycle.
Lemma~\ref{le62} shows that if $Q_3 = Q_1 \circ Q_2$ and $Q_3$ is reduced, then
$Q_3$ is also in the principal cycle.
By Proposition~\ref{pr51} there are integers $k_1 , k_2$ and $k_3$ such that
$\bS_i = \pm \bL_{k_i}$ for $1 \leq i \leq 3$.
What is the relation among the $k_i$'s?  We do not  determine this exactly, but
 show instead the following  approximate additive relation among the $\log ||\bL_{k_i}||$'s.\\

%
%
\begin{lemma}[Infrastructure Bounds]~\label{le63}
 Let $Q_1 , Q_2 , Q_3$ be in the principal cycle and suppose
$\tilde{I}_D \approx Q_1$ via $\pm \bL_{k_1}$,
$\tilde{I}_D \approx Q_2$ via $\pm \bL_{k_2}$, where
$k_1 , k_2 \geq 0$.
Suppose $Q_3 = Q_1 \circ Q_2$ via $\bB$ and that
\beql{eq610}
\log || \bB || \leq c_1 \log D ~.
\eeq
Let $\bS_3$ be defined by
$$
\bS_3 \bB = \bB_0 ( \bS_0 ( \bS_1 \otimes \bS_2 ) ~.
$$ If $\bS_i = \pm \bL_{k_i}$ and $\xi$ is defined by
\beql{eq610a}
  \xi := \log || \bL_{k_1} || + \log || \bL_{k_2} || - \log ||\bL_{k_3} ||,
\eeq
then we have the bound
\beql{eq611}
| \xi | \leq (c_1 + 4) \log D ~. 
\eeq
\end{lemma}

\paragraph{Proof.}
By (\ref{eq608}) we have
\beql{eq612}
|| \bS_3 \bB || = || \bB_0 ( \bS_1 \otimes \bS_2 ) || ~.
\eeq
Now
\begin{eqnarray}
|| \bS_0 ( \bS_1 \otimes \bS_2 ) || & \leq & 4 || \bB_0 || ~|| \bS_1 \otimes \bS_2 || \nonumber \\
~ \nonumber  \\
 & = & 4 || \bB_0 || ~ || \bS_1 || ~ ||\bS_2 || ~.
\end{eqnarray}
We next note that $\bB_0$ is nonnegative and that $\bS_1 \otimes \bS_2$ has
constant sign on columns by Lemma~\ref{le52} (i).
This implies
\beql{eq614}
|| \bB_0 ( \bS_1 \otimes \bS_2 )|| \geq ||\bS_1 \otimes \bS_2 || = || \bS_1 ||~|| \bS_2 || ~.
\eeq
On the other hand
\beql{eq615}
|| \bS_3 \bB || \leq 2 || \bB || ~||\bS_3 ||~.
\eeq
Using orientability the first two columns of $\bB$ form an invertible $2 \times 2$
submatrix $\bB_1$ and we obtain
\beql{eq616}
|| \bS_3 \bB || \geq || \bS_3 \bB_1 || \geq \df{||\bS_3 ||}{2|| \bB_1 ||} \geq
\df{|| \bS_3 ||}{2||\bB_1 ||} \geq \df{||\bS_3 ||}{2||\bB ||} ~,
\eeq
where the center inequality is deduced from
$$
|| \bS_3 || \leq 2 || \bS_3 \bB_1 ||~|| \bB_1^{-1} ||
$$
and
$$
||\bB_1^{-1} || = ( \det B)^{-1} || \bB|| \leq || \bB ||~.
$$
Now (\ref{eq612}), (\ref{eq614}), (\ref{eq615}) yield
\beql{eq617}
|| \bS_3 || \geq \df{2}{|| \bB||} (|| \bS_1 ||~|| \bS_2 ||) ~,
\eeq
while (\ref{eq612}), (6.13), (\ref{eq616}) yield
\beql{eq618}
|| \bS_3 || \leq 8 || \bB_0 ||~|| \bB ||(|| \bS_1 ||~|| \bS_2 ||)~.
\eeq
Using
$$
\log || \bB_0 || \leq \log D
$$
and the hypothesis (\ref{eq610}), the inequalities (\ref{eq617}) and (\ref{eq618})
establish (\ref{eq611}).
$~~\Box$\\

%
%
%

\subsection{Infrastructure Composition Chains for Equivalence on Principal Cycle }

Lemma~\ref{le63} shows that, for $j \ge 1$,  $\log || \bL_j ||$ provides a measure of the
size of the subscript $j$; it shows these quantities are approximately
additive under composition, up to an error (\ref{eq610a}), (\ref{eq611}).
In particular,  composing a form $Q^{(j)}$ with itself
essentially doubles this size.
By repeatedly doubling the size we can rapidly move to forms far apart in the
principal cycle. 
This allows us to find  ``chains" of composition steps going
from $\tilde{I}_D$ to any reduced form in the principal cycle, of length at most
$O((\log D)^2)$ (polynomial in the input size), as given in the following result.
These chains play a role  analogous to "addition chains" in straight-line programming. \\

%
%
\begin{lemma}[Infrastructure Composition Chain] ~\label{le64}
 For any $\bL_j$ with $1 \leq j \leq 2p$ there is a sequence of equivalence matrices
$\bV_k$, and reduced forms $\tilde{Q}_k$ of length $K$ with 
\beql{eq626}
K = O(( \log D)^2 ) ~.
\eeq
having the following properties. 
\begin{description}
\item{(i)}
$\tilde{Q}_0 = \tilde{I}_D , \bV_0 = \left[ \begin{array}{c}
1 ~~ 0 \\
0 ~~ 1
\end{array}
\right]$.
\item{(ii)}
Each pair $( \tilde{Q}_{k+1} , \bV_{k+1} )$ is obtained from the preceding
$( \tilde{Q}_k , \bV_k )$ by a transformation of
either Type I or Type II, where:
\newline
Type I. $\tilde{Q}_{k+1}$ is the right-neighbor of $\tilde{Q}_k$ so that
\beql{eq619}
\tilde \bQ_{k+1} = \bS_{k+1}^T \tilde{Q}_k \bS_{k+1} ~,
\eeq
\beql{eq620}
\bV_{k+1} = \bV_k \bS_{k+1}
\eeq
{\em and}
\beql{eq621}
\log || \bS_k || \leq \df{1}{2} ( \log D ) ~.
\eeq
Type II. $\tilde{Q}_{k+1} = \tilde{Q}_{k_1} \circ \tilde{Q}_{k_2}$ via
$\bB_{k+1}$ for some $0 \leq k_1 , k_2 \leq k$ so that
\beql{eq622}
\bx^T \tilde{\bQ}_i \bx \by^T \tilde{\bQ}_j \by = \bz^T \bB_{k+1}^T \tilde{\bQ}_{k+1} \bB_{k+1}\bz ~,
\eeq
\beql{eq623}
\bV_{k+1} \bB_{k+1} = \bB_0 ( \bV_{k_1} \otimes \bV_{k_2} ) ~,
\eeq
{\em and where}
\beql{eq624}
\log || \bB_{k+1} || = O( \log D ) ~.
\eeq
\item{(iii)}
\beql{eq625}
Q_K = Q^{(j)} \mbox{and $\bV_K = \bL_j$}~.
\eeq
\end{description}
\end{lemma} 

\paragraph{Proof.}
We suppose that the composition of reduced forms is done as in Proposition~\ref{pr61},
so that (\ref{eq624}) is satisfied.
We let $c_1$ denote the constant implied by the O-symbol in (\ref{eq624}).
Let $\sigma_j$ denote the minimal number of type~I and type~II transformations
sequentially applied to get from $\tilde{I}_D $ to $Q^{(j)}$ via $\bL_j$.
First note
\beql{eq627}
\sigma_j \leq j
\eeq
by using type I transformations only.
We will prove by induction on $j$ that for
\beql{eq628}
2(c_1 + 6) \log D \leq j \leq 2p
\eeq
we have
\beql{eq629}
\sigma_j \leq (5+4 (c_1 +4) \log D) ( \log || \bL_j || )~.
\eeq

Suppose (\ref{eq628}) holds.
Take $j_1$ to be some $l$ such that
\beql{eq630}
-(c_1 +5 ) \log D-2 < \log || \bL_l || - \df{1}{2} \log || \bL_j || \leq -(c_1 +4) \log D ~.
\eeq
At least one such $l$ exists by (\ref{eq524})
and $1 \leq l < j$.
(Note (\ref{eq525}) shows $\frac{1}{2} \log || \bL_j || - (c_1+4) \log D \geq 2$.)
Hence we can obtain $\tilde{Q}_k = Q^{(j_1 )} , \bV_k = \bL_{j_1}$ where
$k = \sigma_{j_1}$ satisfies (\ref{eq629}) by the induction hypothesis.
Now apply a type II transformation, using $\tilde{Q}_k \circ \tilde{Q}_k$,
obtaining $\tilde{Q}_{k+1} = Q^{j_2} $ and
$\bV_{k+1} = \bL_{j_2}$.
Using Lemma 6.3 and (\ref{eq630}) we have
$$
2(c_1 + 4) \log D + 2 \geq \log || \bL_j || - \log || \bL_{j_2} || \geq 0 ~.
$$
Then Lemma 5.3 implies
$$
0 \leq j - j_2 < 4 (1+ (c_1 + 4) \log D ) ~.
$$
Hence $4(1+ (c_1 +4) \log D)$ type I transformations will take us to
$Q^{(j)} , \bL _j$.
Hence
\beql{eq631}
\sigma_j \leq \sigma_{j_1} + 4(( c_1 +4) \log D ) +5 ~.
\eeq
But the right side inequality of (\ref{eq630}) gives
\beql{eq632}
\log || \bL_{j_1} || \leq \df{1}{2} \log || \bL_j || \leq 
\log || \bL_j || -1 ~.
\eeq
Substituting (\ref{eq629}) for $j_1$ into (\ref{eq631}) and using (\ref{eq632})
establishes (\ref{eq629}) for $j$ and completes the induction step.
$~~~\Box$ \\
\newline
{\bf Remark.}
By  more detailed  argument, the bound (\ref{eq626}) can be sharpened to
\beql{eq633}
K= O ( \log D) ~.
\eeq
%
%
%

\section{Certificates for Equivalence of Two Indefinite Binary Quadratic Forms}
\hsp
Lemma~\ref{le64} can immediately be used to provide certificates for the equivalence
of two indefinite binary quadratic forms. \\

%
%
\begin{theorem}~\label{th71}
 Let $Q_1$ and $Q_2$ be two indefinite integer binary quadratic forms with 
the same
discriminant.
If $Q_1$ is properly equivalent to $Q_2$, then there is a 
certificate of this equivalence requiring at most
\beql{eq701}
O ( \log || Q_1 || + \log || Q_2 || + ( \log D)^2 M( \log D))
\eeq
elementary operations to verify. 
\end{theorem}

\paragraph{Proof.}
A necessary condition for the equivalence of two forms 
$Q_1 = [a_1 , 2b_1, c_1]$ and
$Q_2 = [a_2 , 2b_2 , c_2 ]$ is that
$$
G.C.D.(a_1 , b_1 , c_1 ) = G.C.D.(a_2 , b_2 , c_2 ) = \sigma_1
$$
and
$$
G.C.D. (a_1 , 2b_1 , c_2 ) = G.C.D.(a_2 , 2b_2 , c_2 ) = \sigma_2 ~.
$$
By removing $\sigma_1$ from the coefficients of both $Q_1$ and $Q_2$ we need only consider
the case $\sigma_1 = 1$.
In that case the forms are {\em properly primitive} if
$\sigma_2 = 1$ and {\em improperly primitive} if $\sigma_2 = 2$.\\

Suppose first that the forms are properly primitive.
Replace $Q_2 = [a_2 , 2b_2 , c_2 ]$ by $\bar{Q}_2 = [a_2 -2b_2 , c_2 ]$, its
{\em inverse} form.
Reduce $Q_1$ and $\bar{Q}_2$, obtaining $Q_1^\ast , \bar{Q}_2^\ast$.
This requires $O( \log || Q_1 || + \log ||Q_2 ||)$ operations by Proposition~\ref{pr42}.
Compose $Q_1^\ast$ and $\bar{Q}_2^\ast$ to obtain a reduced form $Q_3^\ast$.
By Proposition~\ref{pr61} this can be done in $O(M( \log D)) $ operations.\\

Now $Q_1 \sim Q_2$, if and only if $Q_3^\ast \sim \tilde{I}_D$.
This follows from the well-known facts that: (i)  composition of forms induces the
structure of an abelian group on equivalence classes $[Q]$ of properly primitive
forms $Q$, that  (ii) $[ \tilde{I}]_D$ is the identity element of this group, and that (iii) 
$[Q]^{-1} = [ \bar{Q} ]$.
(e.g. see Mathews \cite[Arts. 141, 145]{Mat61}.)\\

We now take the sequence of reduced forms $\tilde{Q}_k$ showing $Q_3^\ast \sim \tilde{I}_D$ that are guaranteed to exist by Lemma~\ref{le64}, together with the
matrices $\bS_k$ and $\bB_k$ involved in the corresponding type~I or II
transformation.
For each transformation we verify either (\ref{eq619}) or (\ref{eq622}), and
this requires $O(M( \log D))$ elementary operations.
We obtain a total of $O(( \log D)^2 M( \log D))$
elementary operations in all, by (\ref{eq626}).\\

Finally, we verify by induction on $k$ that checking (\ref{eq619}), (\ref{eq622}) at each
step guarantees that all $\tilde{Q}_k = \tilde{I}_D$.
Certainly $\tilde{Q}_0 \sim \tilde{I}_D$.
If a type~I transformation is used, then $\tilde{Q}_{k+1} \sim \tilde{Q}_k \sim \tilde{I}_D$
by definition of equivalence.
If a type~II transformation is used, then $\tilde{Q}_i \sim \tilde{I}_D$ and
$\tilde{Q}_j \sim \tilde{I}_D$ guarantees
$\tilde{Q}_{k+1} = \tilde{Q}_1 \circ \tilde{Q}_j \sim \tilde{I}_D$ by Lemma~\ref{le62}.
This completes the proof in the properly primitive case.\\

We treat the improperly primitive case by reducing it to the properly primitive case
by the following method given in 
Mathews \cite[Art. 153]{Mat61}.
We first note that improperly primitive forms have $D \equiv 1 ( \bmod ~4)$.
Let
$$
Q \left| \left[ \begin{array}{c}
a~~b \\
c~~d
\end{array}
\right] \equiv Q(ax+by, cx+dy ) \right.~.
$$
If $D \equiv 1 ( \bmod ~8)$ and $Q$ is improperly primitive, then
$Q \left| \left[ \begin{array}{c}
2 ~~0 \\
0 ~~ 1
\end{array}
\right] = 2Q^\ast \right.$ where $Q$ is properly primitive.
Furthermore if $Q_1 , Q_2$ are two such improperly primitive forms then
$Q_1 \approx Q_2$ if and only if $Q_1^\ast \sim Q_2^\ast$.
We may find a certificate for this as above.
If $D \equiv 5 ( \bmod ~8)$ and $Q$ is improperly primitive,
then $Q \left| \left[ \begin{array}{c}
2 ~ 0 \\
0 ~ 1
\end{array}
\right] = 2Q^{(1)} \right. ,~ Q \left| \left[ \begin{array}{c}
0 ~ 1 \\
0 ~ 2
\end{array}
\right] = 2Q^{(2)}\right. $ and
$Q \left| \left[ \begin{array}{c}
1 ~ 1 \\
0 ~ 2
\end{array}
\right] = 2Q^{(i)}\right. $ where the $Q^{(i)}$ are all properly primitive.
Furthermore if $Q_1 , Q_2$ are two such improperly primitive forms then
$Q_1 \sim Q_2$ if and only if one of $Q_1^{(i)} \sim Q_2^{(1)}$ for $1 \leq i \leq 3$.
We may find a certificate for this as above.
In order to get the bound (\ref{eq701}) we first reduce the improperly primitive forms
and then apply the procedure above. This
reduction uses only $O( \log D)$ additional operations. $~~~\Box$ 

\paragraph{Remark.}  Since $||Q|| > \frac{1}{2} \sqrt D$ for any form $Q$, (\ref{eq701})
gives a bound polynomial in the length of the input $\log || Q_1 || + \log || Q_2 ||$.

%
%
%

\section{Succinct Certificates for BQDE's}
\hsp
We now prove the main results, Theorem~\ref{th11} and Theorem~\ref{th12}.

\paragraph{Proof of Theorem 1.1.}
If (\ref{eq201}) has an admissible solution with $|| \bx || < 256 ||F||^8$, then it
serves as the certificate, and only $O(M( \log |F|))$ operations are needed to
verify it is one.
By Lemma~\ref{le32} this is always the case for definite or degenerate binary
quadratic Diophantine equations.

Now suppose (\ref{eq201}) is indefinite, and has admissible solutions, 
but none
with $|| \bx || < 256 ||F||^8$.
Then by part (ii) of Lemma~\ref{le33}, there exists $\beta_1 , \beta_2$ such that the
\beql{eq801}
y_1^2 - Dy_2^2 = g ~.
\eeq
has a solution with
\beql{eq802}
y_i \equiv \beta_i ( \bmod ~cD \Gamma )
\eeq
satisfying (\ref{eq316}), (\ref{eq317}), and $y_1 > 0$ and the sign of $y_2$
is specified and satisfies one of (\ref{eq315a}), (\ref{eq315b}).
Call the system (\ref{eq801}), (\ref{eq802}) with the given sign conditions 
$E$,
and observe that
\beql{eq803}
||E|| \leq MAX (D, |g|, |cD \Gamma|) \leq 6 || F ||^4
\eeq
using (\ref{eq304}).
By Lemma~\ref{le33} it now suffices to give a certificate for this equation, to
guarantee (\ref{eq300}) has an admissible solution.
Note that it takes only $O(M( \log D) \log D)$ operations to check the
conditions of (ii) of Lemma~\ref{le33}  hold, in particular $O(M ( \log D) \log D)$
operations to compute $\sqrt D$ to one digit past the decimal point, for application to test
the inequalities 
(\ref{eq315a}) and (\ref{eq315b}).
Note $\log ||E|| \le 9 ( \log ||F||)$.\\

Lemma~\ref{le43} shows that to show the system $E$ has an admissible solution it
suffices to produce certificates showing there exist integers $h, B, C$ and
$2 \times 2 $ matrices $\bS , \bW$ such that (i)--(v) of that
lemma hold.
The rest of the proof will accomplish this.\\

Lemma~\ref{le43} also shows that there exist integers $h,B,C$ and $2 \times 2$
matrices $\bS , \bW$ such that (i)--(ix) of that lemma hold.
In the rest of the proof we shall fix this particular choice of $h,B,C, \bS$,
and $\bW$, as well as
\beql{eq804}
Q_{\rm red} = [a_0 , 2b_0 , c_0 ]
\eeq
arising in (iii) of that lemma.
In that case (i), (ii) of Lemma~\ref{le43} can be verified in $O(M( \log ||E||)$
operations by (vi) of that Lemma and (\ref{eq803}).
To verify (iii) of Lemma~\ref{le43} we note that it asserts that
\beql{eq805}
\left[ \begin{array}{c}
a_0 ~~ b_0 \\
b_0 ~~ c_0
\end{array}
\right] = \bS^T \left[
\begin{array}{c}
G~~B \\
B~~C
\end{array}
\right] \bS ~.
\eeq
Using the bound (\ref{eq426}) for a reduced form, (vii), (viii) of Lemma~\ref{le43},
and (\ref{eq803}), all entries in (\ref{eq805}) are
$O( \log ||E||)$ so (\ref{eq805}) can be verified in $O(M( \log ||E||))$ 
operations.\\

The essential difficulty in producing the certificates is the possible large
size of the entries of $\bW$, evidenced by the bound (\ref{eq441}), so that we cannot
afford to keep track of these entries as fixed point binary integers.
Consequently (iv) and (v) of Lemma~\ref{le43} must be verified indirectly.\\

The certificates verifying (iv) and (v) are based on two kinds of formulae, which
we call {\em short} and {\em long}.
The short formulae can be evaluated using fixed-point integer arithmetic.
We will use these to verify (iv).
The long formulae involve integers with too many binary digits to allow
direct evaluation.
We use these to verify (v), by evaluating them $( \bmod ~cD \Gamma )$ to verify
(\ref{eq436}), and by evaluating them using floating-point arithmetic to enough
accuracy to verify (\ref{eq437}).\\

The formulae are those guaranteed to exist by Lemma~\ref{le55} and Lemma \ref{le64}.
By Lemma \ref{le55} the $\bW$ of Lemma \ref{le43} can be written in the form
\beql{eq806}
\bW = (-1)^m ( \bL_{2p} )^K \bL_j
\eeq
for some $j$ with $1 \leq j \leq 2p$, for some $m = 0$ or 1, and  $K$ is bounded by
\begin{eqnarray}~\label{eq807}
|K|&= & O \left(||E||^{\frac{3}{2}}( \log ||E||)^2\right)  \nonumber\\
&=& O\left( ||F||^{17} (\log ||F||)^2\right).
\end{eqnarray}
Assuming that $L_{2p}$ is 
known, we obtain $(\bL_{2p})^K$ by an exponential addition chain of $O\left( \log ||F||\right)$
squarings and multiplications of powers of $\bL_{2p}.$ Then we obtain $\bW$ by
combining this  with
$\bL_j$ using (\ref{eq806}). Here (\ref{eq806}) and the exponential addition
chain formulas are all long formulas.\\

Next, by Lemma \ref{le64} for each $\bL_j$ there exists a chain of
reduced forms $\{ \tilde{Q}_k: ~1 \le k \le K_j\}$ with corresponding reduction matrices $\bS_k$
and  equivalence matrices $\bV_k$
having the properties (\ref{eq619})--(\ref{eq627}). Recall that the type I and II reduction and
composition formulas 
$$
\tilde \bQ_{j+1} = \bS_{j+1}^T\tilde{Q}_j \bS_{j+1}
$$
and
$$
\bx^T \tilde{\bQ}_i \bx \by^T\tilde{\bQ}_j \by = \bz^T \bB_{k+1}^T \tilde{\bQ}_{k+1} \bB_{k+1}\bz
$$
are short formulas,
while the type I and II update formulas 
$$
\bV_{k+1} = \bV_k \bS_{k+1}
$$
and
$$
\bV_{k+1} \bB_{k+1} = \bB_0 ( \bV_{k_1} \otimes \bV_{k_2} ) ~
$$
are long formulas.  \\


Consider the short formulas used in computing $\bL_j$ and $\bL_{2p}$.
Lemma ~\ref{le64}
 gives that all entries in $\bB_k$ and $\bS_k$ have $O(\log D)$ binary
digits.
The size bounds 
(\ref{eq426}) on $\bQ_k$ with these
bounds  imply
that each formula can be evaluated exactly using fixed-point integer arithmetic with
$O ( \log ||E||)$ binary digits.
Each evaluation takes $O(M ( \log ||E||)$ operations, so (\ref{eq626})
implies a total of at most $O(M ( \log ||E||)(\log ||E||)^2)$ bit operations used in
evaluating all the short formulas.
In addition we must verify that the bilinear matrices $\bB_i$ used in short
formulae are unimodular and oriented.
Using the Euclidean algorithm to check unimodularity takes
$O(M( \log ||E||) \log ||E||)$ operations for each $\bB_k$, by \cite[Prop. 3.3]{Lag80a},
for a total of $O(M( \log ||E||)( \log ||E||)^3 )$
operations in all.
Checking orientability requires $O ( M( \log ||E||) \log ||E||)$ 
operations in all.\\

We now verify that the certificate satisfies property (iv) of Lemma~\ref{le43}.
Since $\tilde{I}_D \sim \tilde{I}_D $ via $\bL_{2p}$ and $\tilde{I}_D \sim Q^{(j)}$ via $\bL_j$,
(\ref{eq806})
 implies that
$$
\tilde{I}_D \sim Q^{(j)} ~\mbox{via $\bW$}~.
$$
In order to verify (iv) it suffices to check that
\beql{eq808}
Q^{(j)} = Q_{\rm red},
\eeq
where $Q_{\rm red}$ is as in (\ref{eq433}),
and $Q^{(j)}$ denotes the $Q_k$ produced in (\ref{eq625}) for $\bL_j$.
Checking that (\ref{eq808}) holds takes another $O( \log ||E||)$ operations.

We now describe certificates for (v) of Lemma~\ref{le43}.
We first must verify
\beql{eq809}
hu_i \equiv \alpha_i ~( \bmod ~ cD \Gamma ) ~~ i = 1,2
\eeq
where
\beql{eq810}
\left[ \begin{array}{c}
u_1 ~~ u_3 \\
u_2 ~~ u_4
\end{array}
\right] = \left[
\begin{array}{lr}
1 & - \lambda \\
0 & 1 
\end{array}
\right] \bW \bS^{-1} ~.
\eeq
We define $\bW$ to be given by (\ref{eq806}), and the 
$\bL_{2p}, \bL_j$ are defined
by the long formulae of Lemma~\ref{le64}.
We evaluate all these long formulae as congruences $(\bmod ~cD \Gamma )$.
Since
$$
\log cD \Gamma = O ( \log ||E||) ~,
$$
so we can use binary numbers with $O( \log ||E||)$ digits throughout.
The long formulae for $\bV_k$ in Lemma~\ref{le64} are evaluated successively.
Evaluating each formula (\ref{eq620}) $( \bmod ~cD \Gamma )$ takes
$O( M( \log ||E||)$ operations.
We next must check in (\ref{eq623}) that given $\bB_0 , \bB_{k+1} , \bV_i$
and $\bV_j$ $( \bmod ~cD \Gamma )$ we can calculate 
$\bv_{k+1} ( \bmod ~cD \Gamma )$.
It is straightforward to calculate $\bB_0 ( \bV_{k_i} \otimes \bV_{k_2} )$.
We use the unimodularity condition of the matrix $\bB_{k+1}$, that the greatest
common divisor of its $2 \times 2$ submatrices
$\Delta_{ij} = \left[
\begin{array}{c}
b_{1i} ~~ b_{1j} \\
b_{2i} ~~ b_{2j}
\end{array}
\right]$ is 1.
By an algorithm similar to step~1  of Lagarias \cite[Theorem 5.4]{Lag80a},
repeatedly using the Euclidean algorithm with the $\det ( \Delta_{ij})$
we can find a factorization
\beql{eq811}
c D \Gamma = m_{12} m_{13} m_{14} m_{23} m_{24} m_{34}
\eeq
with the $m_{ij}$ pairwise relatively prime and with
\beql{eq812}
(m_{ij} , \det ( \Delta_{ij} )) = 1 ~.
\eeq
for all $i,j$.
This takes $O(M( \log ||E||) \log ||E||)$ operations.
(Alternatively we can guess a set of $m_{ij}$ and check that they have the
required properties.) Then
\beql{eq813}
( \Delta_{ij} )^{-1} \equiv ( \det ( \Delta_{ij}))^{-1} \left[
\begin{array}{rr}
b_{2j} & -b_{1j} \\
-b_{2i} &  b_{1i}
\end{array}
\right] ( \bmod ~m_{ij} )
\eeq
and $( \det \Delta_{ij} )^{-1} ( \bmod ~m_{ij})$ is calculated in
$O(M( \log ||E||) \log ||E||)$ operations using \cite[Corollary 3.4]{Lag80a}.
Hence
\beql{eq814}
\bV_{k+1} \equiv ( \Delta_{ij} )^{-1} [ \bB_0 ( \bV_{k_1} \otimes 
\bV_{k_2} ) ]_{ij} ( \bmod ~m_{ij} )~,
\eeq
where $[\bM]_{ij}$ denotes the submatrix obtained taking columns $i$ and $j$,
yields $\bV_{k+1} ( \bmod ~m_{ij} )$.
Finally we use the Chinese reminder theorem on each entry of $\bV_{k+1}$
separately to obtain $\bV_{k+1} ( \bmod~ cD \Gamma )$ in 
$O(M ( \log ||E||) \log ||E||)$
operations, by \cite[Prop. 3.6]{Lag80a}.
Thus we may at last obtain $\bL_{2p} , \bL_j ( \bmod ~cD \Gamma )$ in
$O(M (\log ||E||) (\log ||E||)^3 )$ operations, by (\ref{eq626}).
Next we calculate $(\bL_{2p})^2 , (\bL_{2p})^4$ etc.
by successive squarings and reductions $( \bmod ~cD \Gamma )$, and use the
binary expansion of $k$ to evaluate $W( \bmod~ cD \Gamma )$ using 
formula (\ref{eq806})
in $O(M( \log ||E||)(\log ||E||)^2 )$ operations, noting the bound (\ref{eq807}).
Finally (\ref{eq810}) is evaluated $( \bmod ~cD \Gamma )$ and then (\ref{eq809})
verified in a further $O(M( \log ||E||))$ operations.
Thus the congruence conditions
(\ref{eq436}) are  verified in $O(M( \log ||E||)(\log ||E||)^3)$ elementary operations.\\

Finally we check that the sign conditions (\ref{eq437}) of Lemma~\ref{le43} (v) hold.
These can be
verified by evaluating the long formulae using floating-point arithmetic with
floating-point integers maintaining $c_0 ( \log D)^3 = O (( \log ||E||)^3 )$
binary digits in both the exponent and fraction parts, where $c_0$ is a sufficiently
large absolute constant fixed once and for all as described below.
Basic terminology and error estimates for floating-point computations are
given in Appendix~B.
We say that a normalized floating-point number
$\bar{x} = f 2^e$ with $\frac{1}{2} \leq f < 1$ approximates $x$ to
accuracy $s$
{\em significant figures} if
\beql{eq815}
| \bar{x} -x | < 2^{e-s} ~.
\eeq
(Here $(e,f)$ is the representation of $\bar{x}$ used in the calculation.)
We wish to show $u_1$ and $u_2$ are computed to accuracy at least 1 significant
figure, which permits determination of their signs.
Assuming for the moment this accuracy is proved, it is straightforward to
estimate the total number of elementary operations involved in evaluating all the
long formulae to be
$O(M( \log ||E||)^3) ( \log ||E||)^2 )$ which is
$O(M ( \log ||E||) (\log ||E||)^4 )$.
Note here that in evaluating $V_{k+1}$ by (\ref{eq623}) that we merely pick
an invertible $\Delta_{ij}$, and use
\beql{eq816}
\bV_{k+1} = ( \Delta_{ij})^{-1} [ \bB_0 \bV_{k_1} \otimes \bV_{k_2} ]_{ij}
\eeq
evaluated in floating-point, noting that
\beql{eq817}
\log ( \det ( \Delta_{ij} )) = O ( \log ||E|| )
\eeq
using (\ref{eq624}).\\

It remains to estimate the loss of significant figures during the floating-point
computations.
The sources of loss of accuracy in floating-point computations are roundoff error,
exponent overflow, exponent underflow in multiplication, and loss of accuracy in
addition to two nearly equal numbers of opposite signs (e.g. this includes exponent
underflow during addition as a special case).\\

By using $O( \log ||E||)^3 )$ digits in the exponent part, we guarantee that
exponent overflow never occurs.
Indeed, only $O( \log ||E||)$ binary digits are needed to represent the
exponent part $e$ of any entry of $\bW$, since
\beql{eq818}
e = O ( ||E||^{3/2} ( \log ||E||)^2 )
\eeq
by (4.41).
It is easy to check that the bound (\ref{eq818}) applies to any exponent of
every element occurring in the long formulae, since the $V_j$'s
are just various $\bL_k$ with $1 \leq k \leq p$, to which the
bounds (\ref{eq524}), (\ref{eq526}) apply.
Now as long as the floating point calculations agree with the two
entries of the long formulae to one significant figure, their exponents must
agree within $\pm1$ and these calculated exponents will 
then satisfy (\ref{eq818})
and exponent overflow cannot occur.
This demonstrates that exponent overflow cannot occur unless all significant
digits have first been lost due to the other three sources of error.\\

We next show that exponent underflow during multiplications can never occur
 unless all
significant digits have first been lost due to the remaining two sources 
of error.
Indeed the entries of the matrices $\bV_j$ in Lemma \ref{le64} are known a priori 
to be
nonzero integers by Lemma~\ref{le52}(iii),
except for $\bL_i$ with $|i| < 2$ (and if these occur they may be placed at the
beginning of the computation, which is done in fixed point as explained below).
The entries of $\bW$ are nonzero integers since $\bW = \pm \bL_j$ for some
$|j| \geq 2$.
We may suppose the entries of $\bU$ are nonzero integers, for if some
$u_j = 0$ then since $\bU = \bZ$ satisfies the hypotheses of Lemma~\ref{le41} the 
inequalities (\ref{eq409}), (\ref{eq410}) would imply the $u_j$ are small 
enough
that they could be calculated directly in fixed point as
certificates in $O(M( \log ||E||)$ operations to verify (iv), (v) of Lemma~\ref{le43}.
Since these entries are nonzero integers, the exponents of their floating-point
approximations must be $\geq 0$, and exponent underflow during 
multiplication cannot
occur by Lemma~B-1 in Appendix~B.
(We note that some multiplications by zero may occur, but these are exact
using (B--16), (B--17) of Appendix~B.)\\

We must now bound the effects of roundoff error and that of addition of nearly
equal quantities of opposite signs.
We start with $p = c_0 ( \log ||E||)^3$ significant digits of accuracy.
We first consider the calculation of the $\bV_k$ in Lemma~\ref{le63}.
The entries of $\bS_k , \bB_0 , \bB_{k+1}$ are known to $p$ significant
digits by the bounds (\ref{eq621}), (\ref{eq624}).
We will use Lemma~B--2 to bound roundoff error, and Corollary~B--4 to bound
addition of nearly equal quantities.
Evaluating $\bV_{k+1} $ by the long formula (\ref{eq620}) involves a loss of at
most 5 significant digits by Lemma~B--1, since each entry of $\bV_{k+1}$
uses two floating-point multiplications and one addition, and the quantities
added always have the same sign by Lemma~\ref{le52} (i), (ii).
The crucial step lies in showing that evaluating
$\bV_{k+1}$ by the long formula (\ref{eq623}) (actually by (\ref{eq816}) above)
involves a loss of at most $O( \log D)$ significant digits accuracy.
Indeed $\bV_{k_1 } \otimes \bV_{k_2}$ can be evaluated losing at most 3
significant digits accuracy by Lemma~B--1, as only multiplications are
involved.
Now the bound (\ref{eq617}) applies to show that
\beql{eq819}
|| \bV_k || \geq \df{2}{|| \bB_k ||} (|| \bV_{k_1} || ~|| \bV_{k_2} || )|.
\eeq
hence
\beql{eq820}
\log || \bV_k || \geq \log || \bV_{k_1} || + \log || \bV_{k_2} || 
-c_1 \log D -1
\eeq
using (\ref{eq610}).
But all entries of $\bV_k$ have about the same size by Lemma~\ref{le52} (iii), hence
the nearest floating-point approximations to each entry of $\bV_k$ 
must have exponents $e$ satisfying
\beql{eq821}
e \geq \log || \bV_{k_1} || + \log || \bV_{k_2} || - (c_1 +2) \log D-3~.
\eeq
On the
other hand, each entry of $\bV_{k_1} \otimes \bV_{k_2}$ has exponent
\beql{eq822}
e \leq \log || \bV_{k_1} || + \log || \bV_{k_2} || ~.
\eeq
We now evaluate the entries of (\ref{eq816}) doing all multiplications first,
followed by additions.
The multiplications lose at most 6 significant digits each, and the 
resulting exponents satisfy
\beql{eq823}
e \leq \log || \bV_{k_1} || + \log || \bV_{k_2} || + 2c_1 ( \log D) +3
\eeq
using (\ref{eq822}).
Then the additions producing a given entry of $V_{k+1}$ lose at most
$$
(3c_1 +2) \log D + 11
$$
significant digits accuracy, using Corollary~B-4, using (\ref{eq823}) as an
upper bound on $e$ and (\ref{eq821}) as a lower bound on $e - A$.
Thus at most $(8c_1 +2) \log D + 17$ significant digits are lost in
evaluating $\bV_{k+1}$ using the long formula (\ref{eq623}), and thus
at most $O(( \log ||E||)^3 )$ significant digits are lost in evaluating $\bL_{2p}$ and $\bL_j$ 
using Lemma~\ref{le63}.
Next, we note that the calculation of $( \bL_{2p})^k$ in formula (\ref{eq806})
involves a loss of $O( \log D)$ significant digits, because $O( \log D)$ matrix
multiplications are involved in computing $( \bL_{2p})^2 , ( \bL_{2p})^4$ etc., and the
bounds of Lemma~B--2 apply because all numbers added have the same sign.
Calculating $\bW$ using the long formula (\ref{eq806}) loses
another 5 significant digits; again all quantities added have the same sign.
Finally we evaluate
\beql{eq824}
\bU = \left[ \begin{array}{lr}
1 & - \lambda \\
0 & 1
\end{array}
\right] \bW \bS^{-1} ~.
\eeq
where $\lambda = \lceil \sqrt D \rceil$.
Now
$$
\bW = \left[ \begin{array}{c}
1 ~~ \lambda \\
0 ~~ 1
\end{array}
\right] \bU \bS
$$
so
$$
|| \bW || \leq 2 \sqrt D || \bU || ~ || \bS ||
$$
yields
\beql{eq825}
\log || \bU || \geq \log || \bW || -c_2 \log || \bE || ~.
\eeq
for some absolute constant $c_2$, using (\ref{eq440}).
The exponents $e$ of the individual entries of $\bU$ all satisfy
\beql{eq826}
e \geq \log || \bW || -c_3 \log || \bE ||
\eeq
using the inequality (\ref{eq410}) of Lemma~\ref{le41} (which applies since
$\bU = \bZ )$.
Now evaluate the right side of (\ref{eq824}), doing all multiplications first,
and then additions.
The resulting multiplied quantities all have exponents
\beql{eq827}
e \leq \log || \bW || + (c_1 +2 ) \log || \bE || +3~.
\eeq
Then Corollary~B--4 guarantees we can evaluate $\bU$ with a loss of at
most $(c_1 + c_3 + 2) \log ||E|| + 7$ significant digits accuracy.
We have shown at most $O(( \log ||E||)^3 )$ significant digits accuracy can
be lost due to roundoff and adding of nearly equal quantities of
opposite sign in evaluating $\bU$.
Choosing $c_0$ large enough once and for all, we guarantee preservation of
a positive number of significant digits to the end of the computation, and 
Theorem~\ref{th11} is proved.
$~~~\Box$ \\
%
%
\paragraph{ Proof of Theorem~\ref{th12}.}
This essentially follows from Theorem~\ref{th11}. 
The only additional fact that  needs to  be checked is that the certificates of Theorem~\ref{th11} can be
``guessed'' in polynomial time.
The bounds (vi)--(ix) of Lemma~\ref{le43}, the bounds on 
the  the size of the power $k$ in $\bW= \pm(\bL_{2p})^k \bL_j$ in Lemma \ref{le55} and on
the  $\bS_{k+1} , \bB_{k+1}$ in Lemma~\ref{le64} demonstrate that this
can be done.
$~~~\Box$

\newpage

\renewcommand{\theequation}{A.\arabic{equation}}
\setcounter{equation}{0}
\catcode`\@=11
\renewcommand{\section}{
        \setcounter{equation}{1}
        \@startsection {section}{A}{\z@}{-3.5ex plus -1ex minus
        -.2ex}{2.3ex plus .2ex}{\large\bf}
        }
\catcode`@=12

%
%
%

\begin{center}
{\bf Appendix A.~ Period lengths (mod m) of certain linear recurrences.}
\end{center}

Let $(t_1 , u_1 )$ be the least strictly positive solution to Pell's equation
\beql{eqA1}
X^2 - DY^2 = 1
\eeq
and set
\beql{eqA2}
\epsilon = t_1 + u_1 \sqrt D ~.
\eeq
In this appendix we show the sequences $\{t_k \}, \{u_k \}$ defined by
\beql{eqA3}
( \epsilon )^k = t_k + u_k \sqrt D
\eeq
are periodic $( \bmod ~m)$ and we bound the length of the minimal period
$P(m)$ for which
\begin{eqnarray}
t_{k+P(m)} & \equiv & t_k ~( \bmod ~m) \nonumber \\
~~~ \\
u_{k+P(m)} & \equiv & u_k ~( \bmod ~ m)
\end{eqnarray}
both hold.

The sequences $\{t_k\}, \{u_k \}$ both satisfy the second order linear recurrence.
\beql{eqA5}
w_k = t_1 w_{k-1} - w_{k-2} ~.
\eeq
Periodicity of solutions to this recurrence $(\bmod ~m)$ is closely
related to divisibility of $u_k$ by $m$.
Carmichael \cite{Car13}, \cite{Car29} studied divisibility properties of a
class of sequences which includes $\{t_k \} , \{u_k \}$ as special cases.
Periodicity properties for general linear recurrences were considered by
Engstrom \cite{Eng31}, Ward \cite{War33} and other authors. \\

%
%
{\bf Lemma A--1.}
{\em For each $m \ge 1$, the  period $P(m)$ is finite.} 

\paragraph{Proof.}
Let $\bar{\epsilon} = t_1 - u_1 \sqrt D$ so that
\beql{eqA7}
t_k = \df{1}{2} ( \epsilon^k + \epsilon^{-k} )
\eeq
\beql{A8}
u_k = \df{1}{2 \sqrt D} ( \epsilon^k - \epsilon^{-k} )~.
\eeq
Pell's equation asserts that
\beql{A9}
\epsilon \bar{\epsilon} = 1 ~.
\eeq
Thus $\epsilon , \bar{\epsilon}$ are units in the ring of
integers $\sO_D$ of $Q( \sqrt D)$.
For any ideal $\sA$ in $\theta_D$ let $S( \sA )$ denote the smallest
$k$ such that
\beql{eqA10}
\epsilon^k \equiv \bar{\epsilon}^{-k} \equiv 1 ~( \bmod ~ \sA )
\eeq
over $\sO_D$, such $S( \sigma )$ existing since $\epsilon , \bar{\epsilon}$ are
invertible $( \bmod~ \sA )$.
It's easy to check that
\begin{eqnarray*}
t_{k+R} & \equiv & t_k ~( \bmod ~m) \\
~~~ \\
u_{k+R} & \equiv & u_k ~( \bmod ~m)
\end{eqnarray*}
where $R = S (( 2 \sqrt D m))$.
Hence
\beql{eqA11}
P(m)|S((2 \sqrt D m))
\eeq
exists.
$~~~\Box$ \\
%
%

{\bf Lemma A-2.}
{\em If $(m,n) = 1$ then}
\beql{eqA12}
P(mn) = l.c.m. \{P(m), P(n) \}~.
\eeq

\paragraph{Proof.}
This follows from the definition (A-4) and the Chinese Remainder Theorem.
$~~~\Box$ \\

It thus suffices to calculate $P(p^a )$ for prime powers $p^a$. \\
%
%

{\bf Lemma A-3.}
{\em For all primes $p$ and $a \geq 1$,}
\beql{A13}
P(p^{a+1}) | p P(p^a ) ~.
\eeq

\paragraph{Proof.}
For $R = P(p^a )$ we have
\begin{eqnarray}
t_R & = & 1 + p^a  s_1 \nonumber \\
u_R & = & p^a s_2
\end{eqnarray}
for some $s_1 , s_2$.
Since
$$
t_{pR} + u_{pR} \sqrt D = (t_R + u_R \sqrt D )^p
$$
we have, for $p$ an odd prime, 
\beql{eqA15}
t_{pR} = \sum_{j=0}^{p-1/2}
{p \choose 2j}
(t_R )^{p-2j} (u_R )^{2j} D^j
\eeq
\beql{eqA16}
u_{pR} = \sum_{j=0}^{p-1/2}
{p \choose 2j}
(t_R )^{2j} (u_R )^{p-2j} D^{\frac{p-1}{2} -j}
\eeq
Since $p| {p \choose j}$ for $1 \leq j \leq p-1$, these equations and
(A.14) yield
\begin{eqnarray*}
t_{pR} & \equiv & (t_R )^p \equiv 1~ ( \bmod~p^{a+1} ) \\
~~~ \\
u_{pR} & \equiv & 0 ~( \bmod~p^{a+1} ) ~.
\end{eqnarray*}
and (A.13) follows. For the remaining case $p=2$ we have
\begin{eqnarray*}
t_{2R} &= & t_R^2 + D u_R^2 \equiv 1~(\bmod~2^{a+1})\\
u_{2R} &= & 2 t_R u_R ~~~~~\equiv 1~(\bmod~2^{a+1}),
\end{eqnarray*}
giving (A.13) in this case.
$~~~\Box$ \\

In order to bound $P(p)$, let $( \frac{D}{p} )$ denote the Legendre
symbol. \\

%
%

{\bf Lemma A-4.}
{\em Let $p$ be an odd prime.}
\begin{description}
\item{(i)}
{\em If $( \frac{D}{p} ) = 1$, then}
\beql{A17}
P(p)|p-1
\eeq
\item{(ii)}
{\em If $( \frac{D}{p} ) = -1$, then}
\beql{eqA18}
P(p)|2(p+1) ~.
\eeq
\item{(iii)}
{\em If $p|D$, then}
\beql{eqA19}
P(p)|2p ~.
\eeq
\item{(iv)}
$P(2) = 1$ {\em or} 2.
\end{description}

\paragraph{Proof.}
Suppose $p \nmid 2D$ so $( \frac{D}{p}) = \pm 1$.
Then examination of (\ref{eqA7})--(\ref{eqA10}) shows that
(\ref{eqA11}) can be sharpened to
\beql{eqA20}
P(p)|S( p\sO_D)~.
\eeq
\begin{description}
\item{(i)}
If $( \frac{D}{p} ) = 1$, then $(p)$ factors as $(p) = \sP_1 \sP_2$ the
product of two distinct conjugate prime ideals in $\sO_D$.
Then $\sO_D / \sP_i \cong GF(p)$.
Since $x^{p-1} = 1$ in $GF(p)$ when $x \neq 0$,
we have
$$
\epsilon^{p-1} \equiv \bar{\epsilon}^{p-1} \equiv 1 ( \bmod ~ \sP_i )~.
$$
for $i = 1,2$.
Thus
$$
\epsilon^{p-1} \equiv \bar{\epsilon}^{p-1} \equiv 1 ( \bmod ~ p\sO_D)
$$
so $S(p\sO_D)|p-1$.
Then (\ref{eqA20}) proves (A.17).
\item{(ii)}
If $( \frac{D}{p} ) = -1$, then $(p)\sO_D$ is inert, and $\sO_D /(p) \cong GF(p^2 )$.
Now $x^{p+1} \in GF(p)$ for all $x \in GF(p^2 )$ hence
$$
\epsilon{p+1} \equiv a ( \bmod ~ p\sO_D) ~.
$$
for some $a \in Z$.
(Note $GF(p) \cong \ZZ/p\ZZ \subseteq \sO_D /(p)$.)
Applying the conjugation automorphism, we have
$$
\bar{\epsilon}^{p+1} \equiv a ( \bmod ~  p\sO_D)~.
$$
But $\epsilon \bar{\epsilon} = 1$ hence
$$
a^2 = 1 ( \bmod ~ p\sO_D)~.
$$
Hence
$$
\epsilon^{2(p+1)} \equiv \bar{\epsilon}^{2(p+1)} \equiv 1 ( \bmod ~p\sO_D)
$$
and $S(p\sO_D) |2(p+1)$.
Then (A.20) implies (\ref{eqA18}).
\item{(iii)}
If $p|D$ then
$$
t_2 = t_1^2 - Du_1^2 \equiv t_1^2 \equiv 1 ~( \bmod ~p )
$$
since $t_1^2 = 1 + Du_1^2$.
Then (\ref{eqA15}), (\ref{eqA16}) applied with $R=2$ show
\begin{eqnarray*}
t_{2p} & \equiv & (t_2)^p \equiv 1 ~( \bmod ~p) \\
~~~ \\
u_{2p} & \equiv & 0 ~( \bmod ~p ) ~.
\end{eqnarray*}
Hence $P(p)|2p$.
\item{(iv)}
If $2|D$ then $t_1^2 - Du_1^2 = 1$ shows $t_1^2 \equiv 1 ( \bmod ~4)$
and $2|u$, or $4|D$.
In either case $t_2 \equiv 1 ( \bmod ~4)$,
$u_2 \equiv 0~ ( \bmod ~2)$ and $P(2)|2$.
If $2 \nmid D$ then $(t_1 , u_1 ) \equiv (1,0)$ or
(0,1) (mod 2).
In the first case $P(2) =1$, in the second case, the recurrence
(\ref{eqA15}) shows $P(2) = 2$.
$~~~\Box$ \\
\end{description} 

%
%

{\bf Lemma A-5.}
{\em For any $m \ge 2$,}
$$
P(m) \leq 2m( 1+\log m)
$$

\paragraph{Proof.}
Lemmas A-2 through A-4 imply that if $m = \prod_j p_j^{a_j}$ then 
$$
P(m)~ | ~R(m) := 2 \prod_j \left( p_j^{a_j -1} \left( p_j - \left( \df{D}{p} \right) \right)\right).
$$
Now
$$
R(m) \leq 2m \prod_j ( 1 + \df{1}{p_j} ), 
$$
and
\beql{A21}
\prod_j \left( 1 + \df{1}{p_j} \right) \leq \sum_{j=1}^m \df{1}{j} < 1+ \log m,
\eeq
so the lemma follows.
$~~~\Box$ 

\paragraph{Remark.}
By more careful argument one can obtain the improved bound
$P(m) = O (m \log \log m ).$


\renewcommand{\theequation}{B.\arabic{equation}}
\setcounter{equation}{0}
\catcode`\@=11
\renewcommand{\section}{
        \setcounter{equation}{1}
        \@startsection {section}{B}{\z@}{-3.5ex plus -1ex minus
        -.2ex}{2.3ex plus .2ex}{\large\bf}
        }
\catcode`@=12
%
%
%
\begin{center}
{\bf Appendix B.~ Floating-Point Computations.}
\end{center}

This appendix gives upper bounds on the magnitude of errors accumulated in
floating-point computations.
We use the conventions and notation of Knuth [14, Sect.~4.2],
to which we refer for greater detail.\\

We use {\em normalized floating-point numbers} with {\em base} 2,
{\em excess} 0, {\em with p digits}.
Such a number will be denoted $(e,f )$ where
\beql{eqB1}
(e,f) = f2^e ~.
\eeq
Here $e$ is an integer satisfying
\beql{eqB2}
|e| < N
\eeq
and $f$ is a signed fraction such that $2^p f$ is an integer and
satisfying the {\em normalization condition}
\beql{eqB3}
\df{1}{2} \leq |f| < 1 ~.
\eeq
provided $f \neq 0$.
By convention 0 is (0,0).\\

We introduce a notation to distinguish general real numbers from floating-point
numbers, which are just real numbers satisfying (\ref{eqB1})--(\ref{eqB3}).
To this end we always denote floating-point numbers with a bar, i.e., $\bar{x}$
is a floating-point number (to be thought of as an approximation to the
real number $x$).\\

To define the floating-point operations of addition, subtraction, multiplication and
division, we use the function ``Round to $p$ significant figures'' defined by
\beql{eqB4}
\mbox{Round~ $(x,p)$} = \left\{
\begin{array}{ll}
2^{e-p} \lfloor 2^{p-e} x + \df{1}{2} \rfloor , & 2^{e-1} \leq x \leq 2^e \\
~~~ \\
0, & x = 0 \\
~~~ \\
2^{e-p} \lceil 2^{p-e} x- \df{1}{2} \rceil , & 2^{e-1} \leq -x < 2^e~.
\end{array}
\right.
\eeq
We define {\em floating-point addition}
$\oplus $ by
\beql{eqB5}
\bar{x} \oplus \bar{y} = \left\{
\begin{array}{ll}
0, & | \bar{x} + \bar{y} | < 2^{-N} \\
~~~ \\
\mbox{Round~$(x+y,p)$}, & 2^{-E} \leq | \bar{x} + \bar{y} | < 2^N ~.
\end{array}
\right.
\eeq
{\em Exponent overflow} occurs if $| \bar{x} | \bar{y} | \geq 2^N$ and
$\bar{x} \oplus \bar{y}$ is left undefined.
We define floating-point subtraction of $\bar{x}$ as floating-point addition of
$- \bar{x}$.
We define {\em floating-point multiplication} $\oplus$ by
\beql{eqB6}
\bar{x} \otimes \bar{y} = \left\{
\begin{array}{ll}
0, & | \overline{xy}| < 2^{-N} \\
~~~ \\
\mbox{Round~ $(\overline{xy},p)$}, & \mbox{if~$2^{-E} \leq | 
\overline{xy} | < 2^N$}
\end{array}
\right.
\eeq
{\em Exponent overflow} occurs if $| \overline{xy}| \geq 2^N$ and 
$\bar{x} \otimes \bar{y}$
is left undefined.
Floating-point division $\phi$ is
defined similarly to multiplication, but we will not need it.
Note that these operations are well-defined even when {\em exponent underflow}
occurs.\\

Let $\bar{x}$ be a floating point number approximating a nonzero real number
$x$.
Let
\beql{eqB7}
2^e \leq x < 2^{e+1}~.
\eeq
We say $\bar{x}$ {\em approximates} $x$ to $s$ significant digits if
\beql{eqB8}
| \bar{x} -x| < 2^{e-s-1} ~.
\eeq

There are four sources of loss of significant digits in floating-point operations.
\begin{enumerate}
\item
roundoff error,
\item
exponent overflow,
\item
exponent underflow in multiplication,
\item
addition of two nearly equal quantities of opposite signs (includes exponent
underflow).
\end{enumerate}
We deal with these sources separately.\\

Exponent overflow, and exponent underflow in multiplication are the easiest
to handle, by giving sufficient conditions that they do not occur.
By convention multiplication by zero does not count as exponent underflow. \\

%
%
{\bf Lemma B-1.}
{\em Let $\bar{x} = (e_1 , f_1 ), $ and $\bar{y} = (e_2 , f_2 )$ be two floating-point numbers.
If}
\beql{eqB9}
-N + 2 \leq e_1 + e_2 \leq N-1
\eeq
{\em then $\bar{x} \otimes \bar{y}$ does not involve exponent overflow or underflow.
If}
\beql{eqB10}
MAX (e_1 , e_2 ) \leq N-2
\eeq
{\em then $\bar{x} \oplus \bar{y}$ does not involve exponent overflow.}\\

\paragraph{Proof.} Immediate.
$~~~\Box$ \\


In order to analyze roundoff error, we note that when
\beql{eqB11}
2^{e-1} \leq |x| < 2^e
\eeq
we have the bound
\beql{B12}
| \mbox{Round $ (x,p) -x | < 2^{e-p-1}$} ~.
\eeq \\

%
%
{\bf Lemma~B-2.}
{\em Let $\bar{x}- \bar{y}$ be two floating-point numbers, both
having $s$ significant digits.}
\begin{description}
\item{(i)}
{\em If $\bar{x} , \bar{y}$ have the same sign, then at most 2 significant
digits are lost in computing $\bar{x} \oplus \bar{y}$.
\item{(ii)}
If exponent underflow does not occur, then at most 3 significant digits are
lost in computing $\bar{x} \otimes \bar{y}$.}
\end{description}

\paragraph{Proof.}
(i) Since $\bar{x} , \bar{y}$ have the same sign, underflow cannot occur.
Then
\beql{eqB13}
\bar{x} \oplus \bar{y} = \mbox{Round ~ $( \bar{x} + \bar{y} , p)$} ~.
\eeq
Let $\bar{x}, \bar{y}$ have exponents $e_1 , e_2$.
Then the exponent $e_3$ of $\bar{x} \oplus \bar{y}$ is at
least $MAX (e_1 , e_2 )$.
But
\begin{eqnarray}
| \bar{x} - x| & < & 2^{e_1 -s -1} \nonumber \\
~~~ \nonumber \\
| \bar{y} - y | & < & 2^{e_2 -s-1}~.
\end{eqnarray}
Note $s \leq p$.
Then
\begin{eqnarray}
| \bar{x} \oplus \bar{y} -(x+y) | & \leq & | \bar{x} \oplus \bar{y} - 
( \bar{x} + \bar{y}) |+| \bar{x} - x | + | \bar{y} -y | \nonumber \\
~~~ \nonumber \\
 & \leq & 2^{e_3 -p -1} + 2^{e_2 -s-1} +2^{e_1 -s-1} \leq 2^{e_3 -s+1}
\end{eqnarray}
using (B.12), (\ref{eqB13}).\\

(ii)
Since underflow does not occur, we have
\beql{eqB15}
\bar{x} \otimes \bar{y} = \mbox{Round~$( \overline{xy}, p)$}~.
\eeq
If
$e_4$ is the exponent of $\bar{x} \otimes \bar{y}$, then
$$
e_4 \geq e_1 + e_2 - 1 ~.
$$
Now
$$
| \overline{xy} -xy| \leq | \bar{x} -x |~| \bar{y} | + | \bar{y} -y | ~ | x| \leq 2^{e_1 + e_2 -s}
$$
using (B.14).
Hence
\begin{eqnarray*}
| \bar{x} \otimes \bar{y} - xy | & \leq & | \bar{x} \otimes \bar{y} - \overline{xy} |+| \overline{xy} - xy |  \\
~~~ \nonumber \\
 & \leq & 2^{e_3 -p-1} + 2^{e_1 +e_2 -s} \leq 5~2^{e_3 -s-1} ~,
\end{eqnarray*}
using (B.12), (\ref{eqB15}).
$~~~\Box$ \\

We remark that Lemma~B-2 (i) also holds when $\bar{y} = 0$ and
\beql{eqB16}
| \bar{y} - y | < 2^{e_1 -s-1} ~,
\eeq
and that
\beql{eqB17}
\bar{x} \otimes \bar{y} = xy = 0
\eeq
where $y = \bar{y} = 0$.\\

We next consider the bounds for addition. \\

%
%
{\bf Lemma~B-3.}
{\em Let $\bar{x}_1 , \ldots , \bar{x}_j$ be floating-point numbers such that
all $\bar{x}_i$ have exponents $\leq e$.
Suppose that}
\beql{eqB18}
| \bar{x}_i - x_i | < 2^{e-s-1},~~~1 \le i \le j, 
\eeq
{\em and suppose that $e-s \geq -N$.
Let}
\beql{eqB19}
v_j = x_1 + \ldots + x_j ~,
\eeq
{\em and define $\bar{v}_1 = \bar{x}_1$ and}
\beql{eqB20}
\bar{v}_{i+1} = \bar{v}_i \oplus \bar{x}_{i+1},~~~ 2 \le i \le j-1.
\eeq
{\em Let $j \leq 2^k$ and $k \leq p$.
Then}
\beql{eqB21}
| \bar{v}_j -v_j | < 2^{3+2k+3-s} ~.
\eeq

\paragraph{Proof.}
We have
$$
| \bar{x}_i | \leq 2^e - 2^{e-p} ~,
$$
from which it is easy to establish
\beql{eqB22}
\bar{v}_i \leq i( 2^e -2^{e-p} ) (1+ i2^{-p} ) < 2^{e+k+2} ~.
\eeq
(The term $i2^{-p}$ is a roundoff bound.)
Now we have
\beql{eqB23}
| \bar{v}_i -v_i | \leq | \bar{v}_i - ( \bar{v}_{i-1} + \bar{x}_i ) |+| \bar{x}_i
-x_i |+| \bar{v}_{i-1} - v_{i-1} | ~.
\eeq
If we let $e_i$ be the exponent of $\bar{v}_i$ then (B.22) gives
\beql{eqB24}
e_i \leq e + k + 2 ~.
\eeq
But
\beql{eqB25}
| \bar{v}_i - ( \bar{s}_{i-1} - \bar{x}_i ) | \leq MAX (2^{e_i -p -1} , 2^{-N} )~,
\eeq
the bound $2^{-N}$ occurring in the case of underflow.
Then apply (B-18) and (B-25) to (B-22)
and sum over $i$ to obtain
\beql{eqB26}
| \bar{v}_j -v_j | < \sum_{i=1}^j [ 2^{e_j -p-1} +2^{-N} +2^{e-s-1} ]~.
\eeq
using (B-24) gives
\begin{eqnarray*}
| \bar{v}_j -v_j | & < & 2^{e+2k+1-p} + 2^{_N+k} + 2^{e+k-s-1} \\
~~~ \\
  & < & 2^{e+wk+3-s}~,
\end{eqnarray*}
the desired bound.
$~~~\Box$\\

Lemma~B-3 allows one to show that if one knows ``a priori'' that a sum
$\sum_i x_i$ is not too small with respect to its largest term, then the
loss of significant digits in calculating a floating-point approximation
to this sum cannot be large. \\
%
%
{\bf Corollary~B-4.}
{\em Let $\bar{x}_1 , \ldots , \bar{x}_j$ be floating point numbers
approximating $x_1 , \ldots , x_j$ to $s$ significant digits,
with the largest $| \bar{x}_j |$
having exponent $e \geq -N +s$.
Let}
\begin{eqnarray*}
s_j & = & x_1 + \ldots + x_j ~, \\
~~~ \\
\bar{s}_j & = & \bar{s}_{j-1} \oplus \bar{x}_j , \bar{s}_1 = \bar{x}_1 ~,
\end{eqnarray*}
{\em Suppose $j \leq 4$ and that}
$$
|s_j | \geq 2^{e-A} ~.
$$
{\em Then $\bar{s}_j$ approximates $s_j$ to at least $s-A-8$
significant digits.}
$~~~\Box$
%
%
%

\end{document}